\documentclass[12pt]{article}
\usepackage{latexsym,amssymb,amsmath,ulem}
\usepackage{color}
\usepackage{graphicx}
\usepackage{amsfonts}
\usepackage{latexsym}

\def\tto{\;{\lower 1pt \hbox{$\rightarrow$}}\kern -10pt
\hbox{\raise 2pt \hbox{$\rightarrow$}}\;}

\newtheorem{theorem}{Theorem}[section]
\newtheorem{proposition}{Proposition}[section]
\newtheorem{corollary}{Corollary}[section]
\newtheorem{lemma}{Lemma}[section]
\newtheorem{remark}{Remark}[section]
\newtheorem{example}{Example}[section]
\newtheorem{definition}{Definition}[section]

\numberwithin{equation}{section}

\renewcommand{\theequation}{\thesection.\arabic{equation}}
\normalsize

\makeatletter

\ifx\pdfoutput\relax\let\pdfoutput=\undefined\fi
\newcount\msipdfoutput
\ifx\pdfoutput\undefined
\else
 \ifcase\pdfoutput
 \else
    \ifx\paperwidth\undefined
    \else
      \ifdim\paperheight=0pt\relax
      \else
        \pdfpageheight\paperheight
      \fi
      \ifdim\paperwidth=0pt\relax
      \else
        \pdfpagewidth\paperwidth
      \fi
    \fi
  \fi
\fi

\newcount\@hour\newcount\@minute\chardef\@x10\chardef\@xv60
\def\tcitime{
\def\@time{%
  \@minute\time\@hour\@minute\divide\@hour\@xv
  \ifnum\@hour<\@x 0\fi\the\@hour:%
  \multiply\@hour\@xv\advance\@minute-\@hour
  \ifnum\@minute<\@x 0\fi\the\@minute
  }}%

\def\x@hyperref#1#2#3{%
   \catcode`\~ = 12
   \catcode`\$ = 12
   \catcode`\_ = 12
   \catcode`\# = 12
   \catcode`\& = 12
   \catcode`\% = 12
   \y@hyperref{#1}{#2}{#3}%
}

\def\y@hyperref#1#2#3#4{%
   #2\ref{#4}#3
   \catcode`\~ = 13
   \catcode`\$ = 3
   \catcode`\_ = 8
   \catcode`\# = 6
   \catcode`\& = 4
   \catcode`\% = 14
}

\@ifundefined{hyperref}{\let\hyperref\x@hyperref}{}
\@ifundefined{msihyperref}{\let\msihyperref\x@hyperref}{}

\@ifundefined{qExtProgCall}{\def\qExtProgCall#1#2#3#4#5#6{\relax}}{}
\def\QCTOpt[#1]#2{%
  \def\QCTOptB{#1}
  \def\QCTOptA{#2}
}
\def\QCTNOpt#1{%
  \def\QCTOptA{#1}
  \let\QCTOptB\empty
}
\def\Qct{%
  \@ifnextchar[{%
    \QCTOpt}{\QCTNOpt}
}
\def\QCBOpt[#1]#2{%
  \def\QCBOptB{#1}%
  \def\QCBOptA{#2}%
}
\def\QCBNOpt#1{%
  \def\QCBOptA{#1}%
  \let\QCBOptB\empty
}
\def\Qcb{%
  \@ifnextchar[{%
    \QCBOpt}{\QCBNOpt}%
}
\def\PrepCapArgs{%
  \ifx\QCBOptA\empty
    \ifx\QCTOptA\empty
      {}%
    \else
      \ifx\QCTOptB\empty
        {\QCTOptA}%
      \else
        [\QCTOptB]{\QCTOptA}%
      \fi
    \fi
  \else
    \ifx\QCBOptA\empty
      {}%
    \else
      \ifx\QCBOptB\empty
        {\QCBOptA}%
      \else
        [\QCBOptB]{\QCBOptA}%
      \fi
    \fi
  \fi
}
\newcount\GRAPHICSTYPE
\GRAPHICSTYPE=\z@
\def\GRAPHICSPS#1{%
 \ifcase\GRAPHICSTYPE
   \special{ps: #1}%
 \or
   \special{language "PS", include "#1"}%
 \fi
}%
\def\graffile#1#2#3#4{%
    \bgroup
	   \@inlabelfalse
       \leavevmode
       \@ifundefined{bbl@deactivate}{\def~{\string~}}{\activesoff}%
        \raise -#4 \BOXTHEFRAME{%
           \hbox to #2{\raise #3\hbox to #2{\null #1\hfil}}}%
    \egroup
}%
\def\draftbox#1#2#3#4{%
 \leavevmode\raise -#4 \hbox{%
  \frame{\rlap{\protect\tiny #1}\hbox to #2%
   {\vrule height#3 width\z@ depth\z@\hfil}%
  }%
 }%
}%
\newcount\@msidraft
\@msidraft=\z@
\let\nographics=\@msidraft
\newif\ifwasdraft
\wasdraftfalse

\def\GRAPHIC#1#2#3#4#5{%
   \ifnum\@msidraft=\@ne\draftbox{#2}{#3}{#4}{#5}%
   \else\graffile{#1}{#3}{#4}{#5}%
   \fi
}
\def\addtoLaTeXparams#1{%
    \edef\LaTeXparams{\LaTeXparams #1}}%
\newif\ifBoxFrame \BoxFramefalse
\newif\ifOverFrame \OverFramefalse
\newif\ifUnderFrame \UnderFramefalse

\def\BOXTHEFRAME#1{%
   \hbox{%
      \ifBoxFrame
         \frame{#1}%
      \else
         {#1}%
      \fi
   }%
}

\def\doFRAMEparams#1{\BoxFramefalse\OverFramefalse\UnderFramefalse\readFRAMEparams#1\end}%
\def\readFRAMEparams#1{%
 \ifx#1\end%
  \let\next=\relax
  \else
  \ifx#1i\dispkind=\z@\fi
  \ifx#1d\dispkind=\@ne\fi
  \ifx#1f\dispkind=\tw@\fi
  \ifx#1t\addtoLaTeXparams{t}\fi
  \ifx#1b\addtoLaTeXparams{b}\fi
  \ifx#1p\addtoLaTeXparams{p}\fi
  \ifx#1h\addtoLaTeXparams{h}\fi
  \ifx#1X\BoxFrametrue\fi
  \ifx#1O\OverFrametrue\fi
  \ifx#1U\UnderFrametrue\fi
  \ifx#1w
    \ifnum\@msidraft=1\wasdrafttrue\else\wasdraftfalse\fi
    \@msidraft=\@ne
  \fi
  \let\next=\readFRAMEparams
  \fi
 \next
 }%
\def\IFRAME#1#2#3#4#5#6{%
      \bgroup
      \let\QCTOptA\empty
      \let\QCTOptB\empty
      \let\QCBOptA\empty
      \let\QCBOptB\empty
      #6%
      \parindent=0pt
      \leftskip=0pt
      \rightskip=0pt
      \setbox0=\hbox{\QCBOptA}%
      \@tempdima=#1\relax
      \ifOverFrame
          \typeout{This is not implemented yet}%
          \show\HELP
      \else
         \ifdim\wd0>\@tempdima
            \advance\@tempdima by \@tempdima
            \ifdim\wd0 >\@tempdima
               \setbox1 =\vbox{%
                  \unskip\hbox to \@tempdima{\hfill\GRAPHIC{#5}{#4}{#1}{#2}{#3}\hfill}%
                  \unskip\hbox to \@tempdima{\parbox[b]{\@tempdima}{\QCBOptA}}%
               }%
               \wd1=\@tempdima
            \else
               \textwidth=\wd0
               \setbox1 =\vbox{%
                 \noindent\hbox to \wd0{\hfill\GRAPHIC{#5}{#4}{#1}{#2}{#3}\hfill}\\%
                 \noindent\hbox{\QCBOptA}%
               }%
               \wd1=\wd0
            \fi
         \else
            \ifdim\wd0>0pt
              \hsize=\@tempdima
              \setbox1=\vbox{%
                \unskip\GRAPHIC{#5}{#4}{#1}{#2}{0pt}%
                \break
                \unskip\hbox to \@tempdima{\hfill \QCBOptA\hfill}%
              }%
              \wd1=\@tempdima
           \else
              \hsize=\@tempdima
              \setbox1=\vbox{%
                \unskip\GRAPHIC{#5}{#4}{#1}{#2}{0pt}%
              }%
              \wd1=\@tempdima
           \fi
         \fi
         \@tempdimb=\ht1
         \advance\@tempdimb by -#2
         \advance\@tempdimb by #3
         \leavevmode
         \raise -\@tempdimb \hbox{\box1}%
      \fi
      \egroup%
}%
\def\DFRAME#1#2#3#4#5{%
  \vspace\topsep
  \hfil\break
  \bgroup
     \leftskip\@flushglue
	 \rightskip\@flushglue
	 \parindent\z@
	 \parfillskip\z@skip
     \let\QCTOptA\empty
     \let\QCTOptB\empty
     \let\QCBOptA\empty
     \let\QCBOptB\empty
	 \vbox\bgroup
        \ifOverFrame
           #5\QCTOptA\par
        \fi
        \GRAPHIC{#4}{#3}{#1}{#2}{\z@}%
        \ifUnderFrame
           \break#5\QCBOptA
        \fi
	 \egroup
  \egroup
  \vspace\topsep
  \break
}%
\def\FFRAME#1#2#3#4#5#6#7{%
  \@ifundefined{floatstyle}
    {
     \begin{figure}[#1]%
    }
    {
	 \ifx#1h
      \begin{figure}[H]%
	 \else
      \begin{figure}[#1]%
	 \fi
	}
  \let\QCTOptA\empty
  \let\QCTOptB\empty
  \let\QCBOptA\empty
  \let\QCBOptB\empty
  \ifOverFrame
    #4
    \ifx\QCTOptA\empty
    \else
      \ifx\QCTOptB\empty
        \caption{\QCTOptA}%
      \else
        \caption[\QCTOptB]{\QCTOptA}%
      \fi
    \fi
    \ifUnderFrame\else
      \label{#5}%
    \fi
  \else
    \UnderFrametrue%
  \fi
  \begin{center}\GRAPHIC{#7}{#6}{#2}{#3}{\z@}\end{center}%
  \ifUnderFrame
    #4
    \ifx\QCBOptA\empty
      \caption{}%
    \else
      \ifx\QCBOptB\empty
        \caption{\QCBOptA}%
      \else
        \caption[\QCBOptB]{\QCBOptA}%
      \fi
    \fi
    \label{#5}%
  \fi
  \end{figure}%
 }%
\newcount\dispkind%

\def\makeactives{
  \catcode`\"=\active
  \catcode`\;=\active
  \catcode`\:=\active
  \catcode`\'=\active
  \catcode`\~=\active
}
\bgroup
   \makeactives
   \gdef\activesoff{%
      \def"{\string"}%
      \def;{\string;}%
      \def:{\string:}%
      \def'{\string'}%
      \def~{\string~}%
    }
\egroup

\def\FRAME#1#2#3#4#5#6#7#8{%
 \bgroup
 \ifnum\@msidraft=\@ne
   \wasdrafttrue
 \else
   \wasdraftfalse%
 \fi
 \def\LaTeXparams{}%
 \dispkind=\z@
 \def\LaTeXparams{}%
 \doFRAMEparams{#1}%
 \ifnum\dispkind=\z@\IFRAME{#2}{#3}{#4}{#7}{#8}{#5}\else
  \ifnum\dispkind=\@ne\DFRAME{#2}{#3}{#7}{#8}{#5}\else
   \ifnum\dispkind=\tw@
    \edef\@tempa{\noexpand\FFRAME{\LaTeXparams}}%
    \@tempa{#2}{#3}{#5}{#6}{#7}{#8}%
    \fi
   \fi
  \fi
  \ifwasdraft\@msidraft=1\else\@msidraft=0\fi{}%
  \egroup
 }%
\def\TEXUX#1{"texux"}

\def\limfunc#1{\mathop{\rm #1}}%
\def\func#1{\mathop{\rm #1}\nolimits}%
\long\def\QQQ#1#2{%
     \long\expandafter\def\csname#1\endcsname{#2}}%
\@ifundefined{QTP}{\def\QTP#1{}}{}
\@ifundefined{QEXCLUDE}{\def\QEXCLUDE#1{}}{}
\@ifundefined{Qlb}{}{}
\@ifundefined{Qlt}{}{}
\long\def\QQA#1#2{}%
\def\QTR#1#2{{\csname#1\endcsname {#2}}}%

\def\EXPAND#1[#2]#3{}%
\def\NOEXPAND#1[#2]#3{}%
\def\LaTeXparent#1{}%
\def\ChildStyles#1{}%
\def\ChildDefaults#1{}%
\def\QTagDef#1#2#3{}%

\@ifundefined{correctchoice}{}{}
\@ifundefined{HTML}{\def\HTML#1{\relax}}{}
\@ifundefined{TCIIcon}{\def\TCIIcon#1#2#3#4{\relax}}{}
\if@compatibility
\else
  \providecommand{\UNICODE}[2][]{\protect\rule{.1in}{.1in}}
  \providecommand{\U}[1]{\protect\rule{.1in}{.1in}}
  
\fi

\@ifundefined{lambdabar}{
      
   }{}

\@ifundefined{StyleEditBeginDoc}{}{}
\def\QQfnmark#1{\footnotemark}

\@ifundefined{TCIMAKEINDEX}{}{\makeindex}%
\@ifundefined{abstract}{%
 \def\abstract{%
  \if@twocolumn
   \section*{Abstract (Not appropriate in this style!)}%
   \else \small
   \begin{center}{\bf Abstract\vspace{-.5em}\vspace{\z@}}\end{center}%
   \quotation
   \fi
  }%
 }{%
 }%
\@ifundefined{endabstract}{\def\endabstract
  {\if@twocolumn\else\endquotation\fi}}{}%
\@ifundefined{maketitle}{\def\maketitle#1{}}{}%
\@ifundefined{affiliation}{\def\affiliation#1{}}{}%
\@ifundefined{proof}{}{}%
\@ifundefined{endproof}{}{}%
\@ifundefined{newfield}{\def\newfield#1#2{}}{}%
\@ifundefined{chapter}{\def\chapter#1{\par(Chapter head:)#1\par }%
 \newcount\c@chapter}{}%
\@ifundefined{part}{\def\part#1{\par(Part head:)#1\par }}{}%
\@ifundefined{section}{\def\section#1{\par(Section head:)#1\par }}{}%
\@ifundefined{subsection}{\def\subsection#1%
 {\par(Subsection head:)#1\par }}{}%
\@ifundefined{subsubsection}{\def\subsubsection#1%
 {\par(Subsubsection head:)#1\par }}{}%
\@ifundefined{paragraph}{\def\paragraph#1%
 {\par(Subsubsubsection head:)#1\par }}{}%
\@ifundefined{subparagraph}{\def\subparagraph#1%
 {\par(Subsubsubsubsection head:)#1\par }}{}%
\@ifundefined{therefore}{}{}%
\@ifundefined{backepsilon}{}{}%
\@ifundefined{yen}{}{}%
\@ifundefined{registered}{%
   \def\registered{\relax\ifmmode{}\r@gistered
                    \else$\m@th\r@gistered$\fi}%
 \def\r@gistered{^{\ooalign
  {\hfil\raise.07ex\hbox{$\scriptstyle\rm\text{R}$}\hfil\crcr
  \mathhexbox20D}}}}{}%
\@ifundefined{Eth}{}{}%
\@ifundefined{eth}{}{}%
\@ifundefined{Thorn}{}{}%
\@ifundefined{thorn}{}{}%
\@ifundefined{degree}{}{}%
\newdimen\theight
\@ifundefined{Column}{\def\Column{%
 \vadjust{\setbox\z@=\hbox{\scriptsize\quad\quad tcol}%
  \theight=\ht\z@\advance\theight by \dp\z@\advance\theight by \lineskip
  \kern -\theight \vbox to \theight{%
   \rightline{\rlap{\box\z@}}%
   \vss
   }%
  }%
 }}{}%
\@ifundefined{qed}{\def\qed{%
 \ifhmode\unskip\nobreak\fi\ifmmode\ifinner\else\hskip5\p@\fi\fi
 \hbox{\hskip5\p@\vrule width4\p@ height6\p@ depth1.5\p@\hskip\p@}%
 }}{}%
\@ifundefined{cents}{}{}%
\@ifundefined{tciLaplace}{}{}%
\@ifundefined{tciFourier}{}{}%
\@ifundefined{textcurrency}{}{}%
\@ifundefined{texteuro}{}{}%
\@ifundefined{euro}{}{}%
\@ifundefined{textfranc}{}{}%
\@ifundefined{textlira}{}{}%
\@ifundefined{textpeseta}{}{}%
\@ifundefined{miss}{\def\miss{\hbox{\vrule height2\p@ width 2\p@ depth\z@}}}{}%
\@ifundefined{vvert}{}{}
\@ifundefined{tcol}{\def\tcol#1{{\baselineskip=6\p@ \vcenter{#1}} \Column}}{}%
\@ifundefined{dB}{}{}
\@ifundefined{mB}{}{}
\@ifundefined{nB}{}{}
\@ifundefined{note}{}{}%
\def\newfmtname{LaTeX2e}
\ifx\fmtname\newfmtname
  \DeclareOldFontCommand{\rm}{\normalfont\rmfamily}{\mathrm}
  \DeclareOldFontCommand{\sf}{\normalfont\sffamily}{\mathsf}
  \DeclareOldFontCommand{\tt}{\normalfont\ttfamily}{\mathtt}
  \DeclareOldFontCommand{\bf}{\normalfont\bfseries}{\mathbf}
  \DeclareOldFontCommand{\it}{\normalfont\itshape}{\mathit}
  \DeclareOldFontCommand{\sl}{\normalfont\slshape}{\@nomath\sl}
  \DeclareOldFontCommand{\sc}{\normalfont\scshape}{\@nomath\sc}
\fi

\def\alpha{{\Greekmath 010B}}%
\def\beta{{\Greekmath 010C}}%
\def\gamma{{\Greekmath 010D}}%
\def\delta{{\Greekmath 010E}}%
\def\epsilon{{\Greekmath 010F}}%
\def\zeta{{\Greekmath 0110}}%
\def\eta{{\Greekmath 0111}}%
\def\theta{{\Greekmath 0112}}%
\def\iota{{\Greekmath 0113}}%
\def\kappa{{\Greekmath 0114}}%
\def\lambda{{\Greekmath 0115}}%
\def\mu{{\Greekmath 0116}}%
\def\nu{{\Greekmath 0117}}%
\def\xi{{\Greekmath 0118}}%
\def\pi{{\Greekmath 0119}}%
\def\rho{{\Greekmath 011A}}%
\def\sigma{{\Greekmath 011B}}%
\def\tau{{\Greekmath 011C}}%
\def\upsilon{{\Greekmath 011D}}%
\def\phi{{\Greekmath 011E}}%
\def\chi{{\Greekmath 011F}}%
\def\psi{{\Greekmath 0120}}%
\def\omega{{\Greekmath 0121}}%
\def\varepsilon{{\Greekmath 0122}}%
\def\vartheta{{\Greekmath 0123}}%
\def\varpi{{\Greekmath 0124}}%
\def\varrho{{\Greekmath 0125}}%
\def\varsigma{{\Greekmath 0126}}%
\def\varphi{{\Greekmath 0127}}%

\def\nabla{{\Greekmath 0272}}
\def\FindBoldGroup{%
   {\setbox0=\hbox{$\mathbf{x\global\edef\theboldgroup{\the\mathgroup}}$}}%
}

\def\Greekmath#1#2#3#4{%
    \if@compatibility
        \ifnum\mathgroup=\symbold
           \mathchoice{\mbox{\boldmath$\displaystyle\mathchar"#1#2#3#4$}}%
                      {\mbox{\boldmath$\textstyle\mathchar"#1#2#3#4$}}%
                      {\mbox{\boldmath$\scriptstyle\mathchar"#1#2#3#4$}}%
                      {\mbox{\boldmath$\scriptscriptstyle\mathchar"#1#2#3#4$}}%
        \else
           \mathchar"#1#2#3#4%
        \fi
    \else
        \FindBoldGroup
        \ifnum\mathgroup=\theboldgroup 
           \mathchoice{\mbox{\boldmath$\displaystyle\mathchar"#1#2#3#4$}}%
                      {\mbox{\boldmath$\textstyle\mathchar"#1#2#3#4$}}%
                      {\mbox{\boldmath$\scriptstyle\mathchar"#1#2#3#4$}}%
                      {\mbox{\boldmath$\scriptscriptstyle\mathchar"#1#2#3#4$}}%
        \else
           \mathchar"#1#2#3#4%
        \fi     	
	  \fi}

\newif\ifGreekBold  \GreekBoldfalse
\let\SAVEPBF=\pbf
\def\pbf{\GreekBoldtrue\SAVEPBF}%

\@ifundefined{theorem}{\newtheorem{theorem}{Theorem}}{}
\@ifundefined{lemma}{\newtheorem{lemma}[theorem]{Lemma}}{}
\@ifundefined{corollary}{\newtheorem{corollary}[theorem]{Corollary}}{}
\@ifundefined{conjecture}{}{}
\@ifundefined{proposition}{\newtheorem{proposition}[theorem]{Proposition}}{}
\@ifundefined{axiom}{}{}
\@ifundefined{remark}{\newtheorem{remark}{Remark}}{}
\@ifundefined{example}{\newtheorem{example}{Example}}{}
\@ifundefined{exercise}{}{}
\@ifundefined{definition}{\newtheorem{definition}{Definition}}{}

\@ifundefined{mathletters}{%
  \newcounter{equationnumber}
  \def\mathletters{%
     \addtocounter{equation}{1}
     \edef\@currentlabel{\theequation}%
     \setcounter{equationnumber}{\c@equation}
     \setcounter{equation}{0}%
     \edef\theequation{\@currentlabel\noexpand\alph{equation}}%
  }
  
}{}

\@ifundefined{BibTeX}{%
    \def\BibTeX{{\rm B\kern-.05em{\sc i\kern-.025em b}\kern-.08em
                 T\kern-.1667em\lower.7ex\hbox{E}\kern-.125emX}}}{}%
\@ifundefined{AmS}%
    {\def\AmS{{\protect\usefont{OMS}{cmsy}{m}{n}%
                A\kern-.1667em\lower.5ex\hbox{M}\kern-.125emS}}}{}%
\@ifundefined{AmSTeX}{}{}%
\def\@@eqncr{\let\@tempa\relax
    \ifcase\@eqcnt \def\@tempa{& & &}\or \def\@tempa{& &}%
      \else \def\@tempa{&}\fi
     \@tempa
     \if@eqnsw
        \iftag@
           \@taggnum
        \else
           \@eqnnum\stepcounter{equation}%
        \fi
     \fi
     \global\tag@false
     \global\@eqnswtrue
     \global\@eqcnt\z@\cr}

\def\TCItag{\@ifnextchar*{\@TCItagstar}{\@TCItag}}
\def\@TCItag#1{%
    \global\tag@true
    \global\def\@taggnum{(#1)}%
    \global\def\@currentlabel{#1}}
\def\@TCItagstar*#1{%
    \global\tag@true
    \global\def\@taggnum{#1}%
    \global\def\@currentlabel{#1}}
\def\tint{\msi@int\textstyle\int}%
\def\tiint{\msi@int\textstyle\iint}%
\def\tiiint{\msi@int\textstyle\iiint}%
\def\tiiiint{\msi@int\textstyle\iiiint}%
\def\tidotsint{\msi@int\textstyle\idotsint}%
\def\toint{\msi@int\textstyle\oint}%

\def\tsum{\mathop{\textstyle \sum }}%
\def\tbigcup{\mathop{\textstyle \bigcup }}%
\newtoks\temptoksa
\newtoks\temptoksb
\newtoks\temptoksc

\def\msi@int#1#2{%
 \def\@temp{{#1#2\the\temptoksc_{\the\temptoksa}^{\the\temptoksb}}}%
 \futurelet\@nextcs
 \@int
}

\def\@int{%
   \ifx\@nextcs\limits
      \typeout{Found limits}%
      \temptoksc={\limits}%
	  \let\@next\@intgobble%
   \else\ifx\@nextcs\nolimits
      \typeout{Found nolimits}%
      \temptoksc={\nolimits}%
	  \let\@next\@intgobble%
   \else
      \typeout{Did not find limits or no limits}%
      \temptoksc={}%
      \let\@next\msi@limits%
   \fi\fi
   \@next
}%

\def\@intgobble#1{%
   \typeout{arg is #1}%
   \msi@limits
}

\def\msi@limits{%
   \temptoksa={}%
   \temptoksb={}%
   \@ifnextchar_{\@limitsa}{\@limitsb}%
}

\def\@limitsa_#1{%
   \temptoksa={#1}%
   \@ifnextchar^{\@limitsc}{\@temp}%
}

\def\@limitsb{%
   \@ifnextchar^{\@limitsc}{\@temp}%
}

\def\@limitsc^#1{%
   \temptoksb={#1}%
   \@ifnextchar_{\@limitsd}{\@temp}%
}

\def\@limitsd_#1{%
   \temptoksa={#1}%
   \@temp
}

\def\dint{\msi@int\displaystyle\int}%
\def\diint{\msi@int\displaystyle\iint}%
\def\diiint{\msi@int\displaystyle\iiint}%
\def\diiiint{\msi@int\displaystyle\iiiint}%
\def\didotsint{\msi@int\displaystyle\idotsint}%
\def\doint{\msi@int\displaystyle\oint}%

\if@compatibility\else
  \RequirePackage{amsmath}
\fi

\begin{document}

\title{Relaxed Lagrangian duality in convex infinite optimization: reverse
strong duality and optimality}
\author{N. Dinh\thanks{%
International University, VNU-HCM, Linh Trung ward, Thu Duc city, Vietnam
\newline
\null\hskip.5cm (ndinh@hcmiu.edu.vn); and \newline
\null\hskip.5cm Vietnam National University - HCMC, Linh Trung ward, Thu Duc
city, Vietnam} \ \ M. A. Goberna\thanks{%
Department of Mathematics, University of Alicante, Alicante, Spain
(mgoberna@ua.es). Corresponding author.}, \ \ M. A. Lopez\thanks{%
Department of Mathematics, University of Alicante, Alicante, Spain
(marco.antonio@ua.es); and CIAO, Federation University, Ballarat, Australia,
corresponding author}, \ \ M. Volle\thanks{%
Avignon University, LMA EA 2151, Avignon, France
(michel.volle@univ-avignon.fr)} }
\maketitle
\date{}

\begin{abstract}
We associate with each convex optimization problem posed on some locally
convex space with an infinite index set $T,$\ and a given non-empty family $%
\mathcal{H}$ formed by finite subsets of $T,$ a suitable Lagrangian-Haar
dual problem. We provide reverse $\mathcal{H}$-strong duality theorems, $%
\mathcal{H}$-Farkas type lemmas and optimality theorems. Special attention
is addressed to infinite and semi-infinite linear optimization problems.
\end{abstract}


\textit{To Dinh The Luc on the occasion of his 70th anniversary}

\bigskip



\section{Introduction}

In a recent paper on convex infinite optimization \cite{DGLV21}, we have
provided reducibility, zero duality gap, and strong duality theorems for a
new type of Lagrangian-Haar duality associated with families of finite sets
of indices. More precisely, given an optimization problem
\begin{equation}
(\mathrm{P})\ \ \ \ \ \ \inf f(x)\ \ \text{ s.t.}\ \ f_{t}(x)\leq 0,\ \ t\in
T,  \label{1.1}
\end{equation}%
such that $X$ is a locally convex Hausdorff topological vector space, $T$ is
an arbitrary infinite index set, and $\{f;\ f_{t},t\in T\}$ are convex
proper functions on $X$, as well as a family $\mathcal{H}$ of non-empty
finite subsets of the index set $T,$ we consider the $\mathcal{H}$-dual
problem
\begin{equation}
(\mathrm{D}_{\mathcal{H}})\ \ \ \ \ \ \sup\limits_{H\in \mathcal{H},\ \mu
\in \mathbb{R}_{+}^{H}}\inf\limits_{x\in X}\left\{ f(x)+\sum\limits_{t\in
H}\mu _{t}f_{t}(x)\right\} ,  \label{1.3}
\end{equation}%
where $\mu \in \mathbb{R}_{+}^{H}$ stands for $(\mu _{t})_{t\in H}\in
\mathbb{R}_{+}^{H},$\ with the rule $0\times (+\infty )=0.$ When $\mathcal{H}
$ is the family $\mathcal{F}(T)$ of all non-empty finite subsets of $T$, one
gets the standard Lagrangian-Haar dual of $(\mathrm{P}),$
\begin{equation}
(\mathrm{D})\ \ \ \ \ \ \sup\limits_{H\in \mathcal{F}(T),\ \mu \in \mathbb{R}%
_{+}^{H}}\inf\limits_{x\in X}\left\{ f(x)+\sum\limits_{t\in H}\mu
_{t}f_{t}(x)\right\} .  \label{1.2}
\end{equation}%
As in \cite{DGLV21}, this paper pays particular attention to the families $%
\mathcal{H}_{1}:=\left\{ \{t\},\ t\in T\right\} $ of singletons and (when $T=%
\mathbb{N}$) $\mathcal{H}_{\mathbb{N}}:=\left\{ \{1,\ldots ,m\},\ m\in
\mathbb{N}\right\} $ of sets of initial natural numbers. The dual pair $(%
\mathrm{P})-(\mathrm{D}_{\mathcal{H}_{\mathbb{N}}})$ has been used in \cite%
{Karney1983} in the framework of convex semi-infinite programming (CSIP),
where $X=\mathbb{R}^{n}$. More precisely, \cite{Karney1983} gives a
sufficient condition for the optimal value of a SIP problem $(\mathrm{P})$
with $T=\mathbb{N}$ to be the limit, as $m\longrightarrow \infty ,$ of the
optimal values of the sequence of ordinary convex programs $(\mathrm{P}%
_{m})_{m\in \mathbb{N}}$ which results of replacing $T$ by $\{1,\ldots ,m\}$
in $(\mathrm{P}).$ This assumption on $T$ is not as strong as it can seem at
first sight as, if $T$ is an uncountable topological space which contains a
countable dense subset $S$ and the mapping $t\longmapsto f_{t}\left(
x\right) $\ is continuous on $T$ for any $x\in \bigcap_{t\in T}\limfunc{dom}%
f_{t},$ then $(\mathrm{P})$ is equivalent to the countable subproblem which
results of replacing $T$ by $S$ in $(\mathrm{P}).$ In the particular case of
linear semi-infinite programming (LSIP), we can write%
\begin{equation}
(\mathrm{P})\ \ \ \ \ \ \inf \left\langle c^{\ast },x\right\rangle \ \ \text{
s.t.}\ \ \left\langle a_{t}^{\ast },x\right\rangle \leq b_{t},\ \ t\in T,
\label{1.4}
\end{equation}%
with $\left\{ c^{\ast };a_{t}^{\ast },t\in T\right\} \subset \mathbb{R}^{n}$
and $\left\{ b_{t},t\in T\right\} \subset \mathbb{R},$ where, in most
applications, $T$ is a convex body (i.e., a compact convex set with
non-empty interior) in some Euclidean space and the mapping $t\longmapsto
\left( a_{t}^{\ast },b_{t}\right) $\ is continuous on $T$. Then, $T$ can be
replaced by any finite dense subset $S$ to get an equivalent countable LSIP
problem.

There exists a wide literature on the dual pair $(\mathrm{P})-(\mathrm{D}),$
see e.g., the works \cite{DGLS07}, \cite{FLN09}, \cite{FZ16}, \cite{GLV14},
\cite{GLV15}, \cite{LNP08}, and \cite{LZH13}, most of them focused on
constraint qualifications and/or duality theorems, some of them making use,
in order to get optimality conditions, of suitable versions of the
celebrated Farkas' Lemma that have been reviewed in \cite{DJ12}.

The duality theorems for the pair $(\mathrm{P})-(\mathrm{D}_{\mathcal{H}})$
provide conditions guaranteeing a zero duality gap, i.e., that $\inf (%
\mathrm{P})=\sup (\mathrm{D}_{\mathcal{H}})$ (see, \cite[Theorem 6.1]{DGLV21}%
). Other duality theorems in \cite{DGLV21} are strong in the sense that the
optimal value of $(\mathrm{D}_{\mathcal{H}})$ is attained, situation
represented by the equation $\inf (\mathrm{P})=\max (\mathrm{D}_{\mathcal{H}%
})$ (see, \cite[Theorems 5.1-5.3]{DGLV21}). Similarly, the reverse duality
theorems, in Section 3 of this paper, are duality theorems where the optimal
value of $(\mathrm{P})$ is attained, situation represented by the equation $%
\min (\mathrm{P})=\sup (\mathrm{D}_{\mathcal{H}}).$ Reverse (also called
converse) duality theorems for the classical Lagrange dual problem, that is,
for $\mathcal{H}=\mathcal{F}(T),$ in convex infinite programming (CIP in
short) can be found in \cite[Theorem 3.3]{GLV14}{\ and \cite[Theorem 3]%
{GLV15}}. Section 4 provides \textit{ad hoc} Farkas-type results oriented to
obtain, in Section 5, optimality conditions which are expressed in terms of
multipliers associated to the indices belonging to the elements of $\mathcal{%
H}.$

\section{Preliminaries}

Let $X$ be a locally convex Hausdorff topological vector space, and suppose
that its\textbf{\ }topological dual $X^{\ast }$, with null element $%
0_{X^{\ast }},$ is endowed with the weak*-topology. We denote by $\overline{A%
}$ and $\limfunc{ri}A$ the closure and the relative interior of a set $%
A\subset X,$\ and by $\limfunc{co}A$ its convex hull.\ For a set $\emptyset
\neq A\subset X$, by the convex\textbf{\ }cone generated by $A$ we mean $%
\limfunc{cone}A:=\mathbb{R}_{+}(\limfunc{co}A)=\{\mu x:\mu \in \mathbb{R}%
_{+},\ x\in \limfunc{co}A\},$ by $\limfunc{span}A$ its linear span, and by $%
A_{\infty }$ the recession cone of a convex set $A.$ The negative polar of $%
\emptyset \neq A\subset X$ is the convex cone $A^{-}:=\left\{ x^{\ast }\in
X^{\ast }:\left\langle x^{\ast },x\right\rangle \leq 0,\forall x\in
A\right\} .$ The lineality space of a convex cone $K\subset X$ is $\limfunc{%
lin}K=K\cap \left( -K\right) .$

The $w^{\ast }$-closure of a set $\mathbb{A}\subset X^{\ast }$ is also
denoted by $\overline{\mathbb{A}}.$ If $\mathbb{A}\subset X^{\ast }\times
\overline{\mathbb{R}},$ then $\overline{\mathbb{A}}$ denotes the closure of $%
\mathbb{A}$ w.r.t. the product topology. A set $\mathbb{A}\subset X^{\ast
}\times \mathbb{R}$ is said to be $w^{\ast }$-\textit{closed} (respectively,
$w^{\ast }$-closed convex) regarding another subset $\mathbb{B}\subset
X^{\ast }\times \mathbb{R}$\ \ if \ \ $\bar{\mathbb{A}}\cap \mathbb{B}=%
\mathbb{A}\cap \mathbb{B}$ (respectively, $\left( \overline{\limfunc{co}}%
\mathbb{A}\right) \cap \mathbb{B}=\mathbb{A}\cap \mathbb{B}$), see \cite{B10}
(resp. \cite{EV16}).

A function $h:X\rightarrow \overline{\mathbb{R}}:=\mathbb{R}\cup \{\pm
\infty \}$ is proper if its epigraph $\limfunc{epi}h$ is non-empty and never
takes the value $-\infty $; it is convex if $\limfunc{epi}h$ is convex; it
is lower semicontinuous (lsc, in brief) if $\limfunc{epi}h$ is closed; and
it is upper semicontinuous (usc, in brief) if $-h$ is lsc. For a proper
function $h,$ we denote by $[h\leq 0]:=\{x\in X:h(x)\leq 0\}$ its lower
level set of $0,$ and by $\func{dom}h\mathrm{,}$ $\overline{h},$ $\partial
h, $ and $h^{\ast }$ its domain, its lsc envelope,\ its Fenchel
subdifferential, and its Legendre-Fenchel conjugate, respectively. We also
denote by $\Gamma \left( X\right) $ the class of lsc proper convex functions
on $X$. By $\delta _{A}$ we denote the indicator function of $A\subset X,$
with $\delta _{A}\in \Gamma \left( X\right) $ whenever $A\neq \emptyset $ is
closed and convex.

We need to recall some basic facts about convex analysis recession. Given $%
h\in \Gamma (X)$, the recession cone of the closed convex set $\limfunc{epi}%
h $ is the epigraph of the so-called \textit{recession function} $h_{\infty
} $ of $h$: $(\limfunc{epi}h)_{\infty }=\limfunc{epi}h_{\infty }$. The
recession function $h_{\infty }$ coincides with the support function of the
domain of the conjugate $h^{\ast }$ of $h$ (e.g., \cite[Theorem 6.8.5]%
{Laurent-72}):
\begin{equation}
h_{\infty }=\left( \delta _{\func{dom}h^{\ast }}\right) ^{\ast }.
\label{marco1}
\end{equation}%
From \eqref{marco1},
\begin{equation}
\lbrack h_{\infty }\leq 0]=(\func{dom}h^{\ast })^{-}=\{x\in
X\,:\,\left\langle x^{\ast },x\right\rangle \leq 0,\forall x^{\ast }\in
\func{dom}h^{\ast }\},  \label{marco2}
\end{equation}%
which is called the \textit{recession cone} of the function $h$ and provides
the common recession cone to all the non-empty sublevel sets\textbf{\ }$%
[h\leq r]$. Given $\{h_{1},\cdots ,h_{m}\}\subset \Gamma (X)$ such that $%
\bigcap\nolimits_{1\leq k\leq m}\func{dom}h_{k}\neq \emptyset $, by \cite[%
Proposition 3.2.3]{HU-Lema} (whose proof is independent of the dimension of $%
X$), one has for all $\mu \in \mathbb{R}_{+}^{m}:$
\begin{equation}
\left( \sum_{k=1}^{m}\mu _{k}h_{k}\right) _{\infty }=\sum_{k=1}^{m}\mu
_{k}(h_{k})_{\infty }.  \label{marco4}
\end{equation}

\subsection{Classical Lagrange CIP duality}

The \textit{support} of $\lambda :T\rightarrow \mathbb{R}$ is the set $%
\limfunc{supp}\lambda :=\{t\in T:$ $\lambda _{t}\neq 0\}.$ Let $\mathbb{R}%
^{(T)}$ be the \textit{space of generalized finite sequences} formed by all
real-valued functions on $T$ that vanish except on a finite set called
support, i.e.,
\begin{equation*}
\mathbb{R}^{(T)}:=\{\lambda :T\rightarrow \mathbb{R}_{+}\text{ such that }%
\limfunc{supp}\lambda \text{ is finite}\},
\end{equation*}%
\ with positive cone $\mathbb{R}_{+}^{(T)}:=\{\lambda \in \mathbb{R}^{(T)}:$
$\lambda _{t}\geq 0,\forall t\in T\}.$ We can associate to each $\lambda \in
\mathbb{R}_{+}^{(T)}$ the function\textbf{\ }$\sum_{t\in T}\lambda
_{t}f_{t}\colon \ X\rightarrow \mathbb{R}\cup \{+\infty \}$ such that
\begin{equation*}
\left( \sum_{t\in T}\lambda _{t}f_{t}\right) (x)=\left\{
\begin{tabular}{ll}
$\sum\limits_{t\in \limfunc{supp}\lambda }\lambda _{t}f_{t}(x),\text{ }$ & $%
\text{if }\limfunc{supp}\lambda \neq \emptyset ,$ \\
$0,$ & $\text{if }\limfunc{supp}\lambda =\emptyset .$%
\end{tabular}%
\right.
\end{equation*}%
So, we can reformulate $(\mathrm{D})$ in (\ref{1.2})\ as
\begin{equation*}
(\mathrm{D})\quad \sup_{\lambda \in \mathbb{R}_{+}^{(T)}}\inf_{x\in
X}\left\{ f(x)+\left( \sum_{t\in T}\lambda _{t}f_{t}\right) (x)\right\} .
\end{equation*}%
It is known that the function\textbf{\ }$\varphi :X^{\ast }\rightarrow
\overline{\mathbb{R}}$ such that \textbf{\ }%
\begin{equation*}
\varphi (x^{\ast }):=\inf_{\lambda \in \mathbb{R}_{+}^{(T)}}\left(
f+\sum_{t\in T}\lambda _{t}f_{t}\right) ^{\ast }(x^{\ast })
\end{equation*}%
and the set \textbf{\ }%
\begin{equation*}
\mathcal{A}:=\bigcup_{\lambda \in \mathbb{R}_{+}^{(T)}}\limfunc{epi}\left(
f+\sum_{t\in T}\lambda _{t}f_{t}\right) ^{\ast }\subset X^{\ast }\times
\mathbb{R}
\end{equation*}%
are both convex, and $\limfunc{epi}\overline{\varphi }=\overline{\mathcal{A}}
$ (see, for instance, \cite{DGLV21}, \cite{GLV14}, \cite{GLV15})

We denote the feasible set of $(\mathrm{P})$ by
\begin{equation*}
E:=\bigcap_{t\in T}[f_{t}\leq 0].
\end{equation*}%
Then,
\begin{equation*}
-\infty \leq (f+\delta _{E})^{\ast }(x^{\ast })\leq \varphi (x^{\ast })\leq
f^{\ast }(x^{\ast })\leq +\infty ,\;\forall x^{\ast }\in X^{\ast },
\end{equation*}%
and, taking $x^{\ast }=0_{X^{\ast }}$, one gets the weak duality for the
pair $(\mathrm{P})-(\mathrm{D}):$
\begin{equation*}
-\infty \leq \inf\nolimits_{X}f\leq \sup (\mathrm{D})\leq \inf (\mathrm{P}%
)\leq +\infty .
\end{equation*}

\subsection{Relaxed Lagrange CIP duality}

Let $\mathcal{H}$ be a non-empty family of non-empty finite subsets of $T,$
that is, $\emptyset \neq \mathcal{H}\subset \mathcal{F}(T),$ with associated
dual problem $(\mathrm{D}_{\mathcal{H}})$ as in (\ref{1.3}). Obviously,
\begin{equation}
\sup (\mathrm{D}_{\mathcal{H}})\leq \sup (\mathrm{D}_{\mathcal{F}(T)})=\sup (%
\mathrm{D})\leq \inf (\mathrm{P}).  \label{2.5}
\end{equation}%
Let us define the sets%
\begin{eqnarray*}
E_{\mathcal{H}}:= &&\bigcap\limits_{H\in \mathcal{H},t\in H}[f_{t}\leq 0], \\
\mathcal{A}_{\mathcal{H}}:= &&\bigcup\limits_{H\in \mathcal{H},\mu \in
\mathbb{R}_{+}^{H}}\limfunc{epi}\left( f+\sum_{t\in H}\mu _{t}f_{t}\right)
^{\ast },
\end{eqnarray*}%
and the function $\varphi _{\mathcal{H}}:X^{\ast }\longrightarrow \mathbb{R}$
such that
\begin{equation*}
\varphi _{\mathcal{H}}:=\inf\limits_{H\in \mathcal{H},\mu \in \mathbb{R}%
_{+}^{H}}\left( f+\sum_{t\in H}\mu _{t}f_{t}\right) ^{\ast }.
\end{equation*}%
Obviously, $A_{\mathcal{H}}\subset A$ and $\varphi _{\mathcal{H}}\geq
\varphi .$

\begin{definition}
$\mathrm{(i)}$ A family $\mathcal{H}\subset \mathcal{F}(T)$ is said to be
\textit{covering} if $\bigcup\nolimits_{H\in \mathcal{H}}H=T$.\newline
$\mathrm{(ii)}$ A family $\mathcal{H}\subset \mathcal{F}(T)$ is said to be
\textit{directed} if for each $H,K\in \mathcal{H}$ there exists $L\in
\mathcal{H}$ such that $H\cup K\subset L$.
\end{definition}

The families $\mathcal{F}(T)$ and $\mathcal{H}_{\mathbb{N}}$ are both
covering and directed families, whereas\textbf{\ }$\mathcal{H}_{1}$ is just
covering.

As shown in \cite[Proposition 3.2]{DGLV21}, for each directed covering
family $\mathcal{H}\subset \mathcal{F}(T)$ one has
\begin{equation}
\mathcal{A}_{\mathcal{H}}=\mathcal{A}_{\mathcal{F}(T)}=\mathcal{A},
\label{2.1}
\end{equation}%
and, consequently,
\begin{equation}
\varphi _{\mathcal{H}}=\varphi _{\mathcal{F}(T)}=\varphi ,\;\text{\ and }%
\sup (\mathrm{D}_{\mathcal{H}})=\sup (\mathrm{D}_{\mathcal{F}(T)})\equiv
\sup \mathrm{(D}).  \label{2.2}
\end{equation}%
Let $\mathcal{H}\subset \mathcal{F}(T)$ be a covering family. Then, $E_{%
\mathcal{H}}=E$ and, according to \cite[Lemma 5.2]{DGLV21}, $\{f;\
f_{t},t\in T\}\subset \Gamma (X)$ entails
\begin{equation}
(\varphi _{\mathcal{H}})^{\ast }=f+\delta _{E}  \label{2.3}
\end{equation}%
and if, additionnally, $E\cap (\func{dom}f)\neq \emptyset $, then
\begin{equation*}
\limfunc{epi}(f+\delta _{E})^{\ast }=\overline{\limfunc{co}}\mathcal{A}_{%
\mathcal{H}}=\overline{\limfunc{co}}\left( \bigcup\limits_{H\in \mathcal{H}%
,\mu \in \mathbb{R}_{+}^{H}}\limfunc{epi}\left( f+\sum_{t\in H}\mu
_{t}f_{t}\right) ^{\ast }\right) .
\end{equation*}%
Moreover, by \cite[Theorem 5.1]{DGLV21}, $\mathcal{H}$-strong duality holds
at a given $x^{\ast }\in X^{\ast },$ i.e.,%
\begin{equation}
(f+\delta _{E})^{\ast }(x^{\ast })=\min_{H\in \mathcal{H}\text{,\ }\mu \in
\mathbb{R}_{+}^{H}}\left( f+\sum_{t\in H}\mu _{t}f_{t}\right) ^{\ast
}(x^{\ast }),  \label{4.4}
\end{equation}%
\ if and only if $\mathcal{A}_{\mathcal{H}}$ is $w^{\ast }$-closed convex
regarding $\{x^{\ast }\}\times \mathbb{R}.$

\subsection{The $\mathcal{H}$-dual problem as a limit}

It is easy to see that the mapping $\mathcal{F}(T)\supset \mathcal{H}%
\longmapsto \sup (\mathrm{D}_{\mathcal{H}})\in \overline{\mathbb{R}}$ is
non-decreasing w.r.t. the inclusion $\subset $ in $\mathcal{F}(T).$
Consequently, if the family $\mathcal{H}\subset \mathcal{F}(T)$ is directed,
we can express $\sup (\mathrm{D}_{\mathcal{H}})$ as the limit of a net as
follows:%
\begin{equation*}
\sup (\mathrm{D}_{\mathcal{H}})=\sup\limits_{H\in \mathcal{H}}\sup (\mathrm{D%
}_{H})=\lim_{H\in \mathcal{H}}\sup (\mathrm{D}_{H}).
\end{equation*}%
If, moreover, $\mathcal{H}$ is covering, then
\begin{equation}
\sup (\mathrm{D})=\lim_{H\in \mathcal{H}}\sup (\mathrm{D}_{H}).  \label{2.6}
\end{equation}
In particular, if $T=\mathbb{N},$ we consider the countable program
\begin{equation}
(\mathrm{P}_{\mathbb{N}})\quad \inf f(x)\,\,\mathrm{s.t.}\,f_{k}(x)\leq
0,k\in \mathbb{N},  \label{2.7}
\end{equation}%
and the sequence of finite subproblems
\begin{equation}
(\mathrm{P}_{m})\quad \inf f(x)\,\,\mathrm{s.t.}\,f_{k}(x)\leq 0,\ k\in
\{1,\cdots ,m\},\ m\in \mathbb{N},  \label{2.8}
\end{equation}%
whose ordinary Lagrangian dual problems are
\begin{equation}
(\mathrm{D}_{m})\ \ \ \ \ \ \sup\limits_{\mu \in \mathbb{R}%
_{+}^{m}}\inf\limits_{x\in X}\left\{ f(x)+\sum\limits_{k=1}^{m}\mu
_{k}f_{k}(x)\right\} ,\ m\in \mathbb{N}.  \label{2.9}
\end{equation}%
From (\ref{2.6}), the Lagrangian-Haar dual of $(\mathrm{P}_{\mathbb{N}}),$
\begin{equation}
(\mathrm{D}_{\mathbb{N}})\quad \sup_{\lambda \in \mathbb{R}_{+}^{(\mathbb{N}%
)}}\inf_{x\in X}\left\{ f(x)+\sum_{k\in \mathbb{N}}\lambda
_{k}f_{k}(x)\right\} ,  \label{2.10}
\end{equation}%
and its $\mathcal{H}_{\mathbb{N}}$-dual Lagrange problem $(\mathrm{D}_{%
\mathcal{H}_{\mathbb{N}}})$ can be expressed as limits in this way:
\begin{equation}
\sup (\mathrm{D}_{\mathbb{N}})=\sup (\mathrm{D}_{\mathcal{H}_{\mathbb{N}%
}})=\lim_{m\rightarrow \infty }\sup (\mathrm{D}_{m}).  \label{2.11}
\end{equation}%
Corollary \ref{cor7} below provides a sufficient condition for the primal
counterpart of (\ref{2.9}):
\begin{equation*}
\inf (\mathrm{P}_{\mathbb{N}})=\lim_{m\rightarrow \infty }\inf (\mathrm{P}%
_{m}).
\end{equation*}

\section{$\mathcal{H}$-reverse strong duality}

Let us go back to the general convex infinite optimization problem $(\mathrm{%
P})$ in (\ref{1.1}). Along this section we assume that $\{f;\ f_{t},t\in
T\}\subset \Gamma (X)$ and $E\cap \func{dom}f\neq \emptyset $, meaning that $%
\inf (\mathrm{P})\neq +\infty $.

\begin{definition}
Given a covering family $\mathcal{H}\subset \mathcal{F}(T)$, we say that $%
\mathcal{H}$-\textit{reverse strong duality} holds if
\begin{equation*}
\min (\mathrm{P})=\sup (\mathrm{D}_{\mathcal{H}}),
\end{equation*}%
equivalently, that\textbf{\ }there exists $\bar{x}\in E \cap \func{dom} f$
such that
\begin{equation*}
f(\bar{x})=\sup (\mathrm{D}_{\mathcal{H}})\in \mathbb{R}.
\end{equation*}
\end{definition}

We first show that $\mathcal{H}$-reverse strong duality can be described in
terms of subdifferentiability of the function $\varphi _{\mathcal{H}}$.%
\newline

Recall that the subdifferential of a function $g:X^{\ast }\rightarrow
\overline{\mathbb{R}}$ at a point $a^{\ast }\in X^{\ast }$ is given by
\begin{equation*}
\partial g(a^{\ast }):=%
\begin{cases}
\left\{ x\in X\,:\,g(x^{\ast })\geq g(a^{\ast })+\left\langle x^{\ast
}-a^{\ast },x\right\rangle ,\,\forall x^{\ast }\in X^{\ast }\right\} , &
\text{ if }g(a^{\ast })\in \mathbb{R}, \\
\emptyset , & \text{ if }g(a^{\ast })\notin \mathbb{R}.%
\end{cases}%
\end{equation*}%
We have
\begin{equation}
x\in \partial g(a^{\ast })\Leftrightarrow g(a^{\ast })+g^{\ast
}(x)=\left\langle a^{\ast },x\right\rangle .  \label{3.1}
\end{equation}

\begin{lemma}
\label{lem3} Let $\mathcal{H}$ be a covering family. Then, $\mathcal{H}$%
-reverse strong duality holds if and only if $\varphi _{\mathcal{H}}$ is
subdifferentiable at $0_{X^{\ast }}$. In such a case one has $\partial
\varphi _{\mathcal{H}}(0_{X^{\ast }})=\limfunc{sol}(\mathrm{P}),$ where%
\textbf{\ }$\limfunc{sol}(\mathrm{P})$ is\textbf{\ }the optimal solution set
of $(\mathrm{P})$.
\end{lemma}

\noindent \textbf{Proof} Let $x\in \partial \varphi _{\mathcal{H}%
}(0_{X^{\ast }}).$ Since we are assuming that $\mathcal{H}$ is covering,\ by
(\ref{2.3}) and (\ref{3.1}),\textbf{\ }we have
\begin{equation*}
(f+\delta _{E})(x)=(\varphi _{\mathcal{H}})^{\ast }(x)=-\varphi _{\mathcal{H}%
}(0_{X^{\ast }})\in \mathbb{R}.
\end{equation*}%
Then $x\in E$ and
\begin{equation*}
\inf (\mathrm{P})\leq f(x)=-\varphi _{\mathcal{H}}(0_{X^{\ast }})=\sup (%
\mathrm{D}_{\mathcal{H}})\leq \inf (\mathrm{P}).
\end{equation*}%
Consequently, if\textbf{\ }$\varphi _{\mathcal{H}}$ is subdifferentiable at $%
0_{X^{\ast }}$ then $\mathcal{H}$-reverse strong duality holds and $\partial
\varphi _{\mathcal{H}}(0_{X^{\ast }})\subset \limfunc{sol}(\mathrm{P})$.
\newline
Assume now that $\mathcal{H}$-reverse strong duality holds. There exists $%
x\in E\cap (\func{dom}f)$ such that
\begin{equation}
(\varphi _{\mathcal{H}})^{\ast }(x)=f(x)=\sup (\mathrm{D}_{\mathcal{H}%
})=-\varphi _{\mathcal{H}}(0_{X^{\ast }})\in \mathbb{R},  \label{3.2}
\end{equation}%
that means $x\in \partial \varphi _{\mathcal{H}}(0_{X^{\ast }})$ and the
first part of Lemma \ref{lem3} is proved with, in addition, the inclusion $%
\partial \varphi _{\mathcal{H}}(0_{X^{\ast }})\subset \limfunc{sol}(\mathrm{P%
})$. It remains to prove that if $\mathcal{H}$-reverse strong duality holds,
then $\limfunc{sol}(\mathrm{P})\subset \partial \varphi _{\mathcal{H}%
}(0_{X^{\ast }})$. Now for each $x\in \limfunc{sol}(\mathrm{P})$ we have (%
\ref{3.2}). So, $\varphi _{\mathcal{H}}(0_{X^{\ast }})+(\varphi _{\mathcal{H}%
})^{\ast }(x)=0,$ that means $x\in \partial \varphi _{\mathcal{H}%
}(0_{X^{\ast }})$. \hfill $\square \medskip \smallskip $

In favorable circumstances we know that $\varphi _{\mathcal{H}}$ is a convex
function. For instance, when the covering family $\mathcal{H}$ is also
directed, by (\ref{2.1}) and (\ref{2.2}), $\mathcal{A}_{\mathcal{H}}=%
\mathcal{A}$ and $\varphi _{\mathcal{H}}=\varphi ,$ respectively, implying
the convexity of both $\mathcal{A}_{\mathcal{H}}$ and $\varphi _{\mathcal{H}%
}.$ Another important example is furnished by
\begin{equation*}
\varphi _{\mathcal{H}_{1}}=\inf_{(t,\mu )\in T\times \mathbb{R}_{+}}(f+\mu
f_{t})^{\ast },
\end{equation*}%
which is convex under the assumptions \textrm{(a)}, \textrm{(b)}, \textrm{(c)%
} of Corollary \ref{cor5} below (see \cite[Remark 5.5]{DGLV21}). In order to
propose a tractable subdifferentiability criterion when $\varphi _{\mathcal{H%
}}$ is convex we need to recall some facts about quasicontinuous convex
functions and convex analysis recession.\newline

\begin{definition}
A convex function $g:X^{\ast }\rightarrow \overline{\mathbb{R}}$ is said to
be $\tau (X^{\ast },X)$-quasicontinuous (\cite{Zoly1970}, \cite{Zoly1971}),
where $\tau $ is the Mackey topology on $X^{\ast },$ if the following four
properties are satisfied:
\end{definition}

\begin{enumerate}
\item $\mathrm{aff}(\func{dom}g)$ is $\tau (X^{\ast },X)$-closed (or $\sigma
(X^{\ast },X)$-closed),

\item $\mathrm{aff}(\func{dom}g)$ is of finite codimension,

\item the $\tau (X^{\ast },X)$-relative interior of $\func{dom}g$, say $%
\mathrm{ri}(\func{dom}g)$, is non-empty,

\item the restriction of $g$ to $\mathrm{aff}(\func{dom}g)$ is $\tau
(X^{\ast },X)$-continuous on $\mathrm{ri}(\func{dom}g)$.
\end{enumerate}

Lemmas \ref{lem4}, \ref{lem5}, \ref{lem6} below will be used in the sequel.

\begin{lemma}[{\protect\cite[Proposition 5.4]{Zoly1970}}]
\label{lem4} Let $h\in \Gamma (X)$. The conjugate function $h^{\ast }$ is $%
\tau (X^{\ast },X)$-quasicontinuous if and only if $h$ is weakly inf-locally
compact; that is to say $[h\leq r]$ is weakly locally compact for each $r\in
\mathbb{R}$.
\end{lemma}

\begin{lemma}[{\protect\cite[Theorem II.4]{Volle1997}}]
\label{lem5} A convex function $g:X^{\ast }\rightarrow \overline{\mathbb{R}}$
majorized by a $\tau (X^{\ast },X)$-quasicontinuous one is $\tau (X^{\ast
},X)$-quasicontinuous, too.
\end{lemma}

\begin{lemma}[{\protect\cite[Theorem III.3]{Volle1997}}]
\label{lem6} Let $g:X^{\ast }\rightarrow \overline{\mathbb{R}}$ be a $\tau
(X^{\ast },X)$-quasicontinuous convex function such that $g(0_{X^{\ast
}})\neq -\infty $ and $\overline{\limfunc{cone}}\func{dom}g$ is a linear
subspace of $X^{\ast }$. Then $\partial g(0_{X^{\ast }})$ is the sum of a
non-empty weakly compact convex set and a finite dimensional linear subspace
of $X$.
\end{lemma}

We define the {recession cone of }$(\mathrm{P})$ by setting
\begin{equation*}
(\mathrm{P})_{\infty }:=\bigcap_{t\in T}[(f_{t})_{\infty }\leq 0]\cap
\lbrack f_{\infty }\leq 0].
\end{equation*}
For the next theorem and the corollaries below, recall that $\inf \mathrm{(P)%
} \not= + \infty$ as $E\cap \func{dom}f\neq \emptyset $.

\begin{theorem}[$\mathcal{H}$-reverse strong duality]
\label{thm4} Let $\mathcal{H}$ be a covering family such that $\varphi _{%
\mathcal{H}}$ is convex $\tau (X^{\ast },X)$-quasicontinuous and $(\mathrm{P}%
)_{\infty }$ is a linear subspace of $X$. Then $\mathcal{H}$-reverse strong
duality holds:
\begin{equation*}
\min (\mathrm{P})=\sup (\mathrm{D}_{\mathcal{H}})\in \mathbb{R}.
\end{equation*}%
Moreover, $\limfunc{sol}(\mathrm{P})$ is the sum of a weakly compact convex
set and a finite dimensional linear subspace of $X$.
\end{theorem}

\noindent \textbf{Proof }One has $\varphi _{\mathcal{H}}(0_{X^{\ast
}})=-\sup (\mathrm{D}_{\mathcal{H}})\geq -\inf (\mathrm{P})>-\infty $ (the
last strict inequality holds as $E\cap \func{dom}f\neq \emptyset $)$. $ In
order to apply Lemma \ref{lem6} to the convex function $\varphi _{\mathcal{H}%
}$, we have to prove that $\overline{\limfunc{cone}}\func{dom}\varphi _{%
\mathcal{H}}$ is a linear subspace. We have
\begin{equation*}
\overline{\limfunc{cone}}\func{dom}\varphi _{\mathcal{H}}=(\func{dom}\varphi
_{\mathcal{H}})^{--}=\{x^{\ast }\in X^{\ast }\,:\,\left\langle x^{\ast
},x\right\rangle \leq 0,\forall x\in (\func{dom}\varphi _{\mathcal{H}%
})^{-}\}.
\end{equation*}%
Therefore, $\overline{\limfunc{cone}}\func{dom}\varphi _{\mathcal{H}}$ is a
linear subspace if and only if $(\func{dom}\varphi _{\mathcal{H}})^{-}$ is a
linear subspace. Now,
\begin{equation*}
\func{dom}\varphi _{\mathcal{H}}=\bigcup\limits_{H\in \mathcal{H}%
}\bigcup\limits_{\mu \in \mathbb{R}_{+}^{H}}\func{dom}\left( f+\sum_{t\in
H}\mu _{t}f_{t}\right) ^{\ast }
\end{equation*}%
and we can write%
\begin{align*}
(\func{dom}\varphi _{\mathcal{H}})^{-}& =\bigcap_{H\in \mathcal{H}%
}\bigcap_{\mu \in \mathbb{R}_{+}^{H}}\left( \func{dom}\left( f+\sum_{t\in
H}\mu _{t}f_{t}\right) ^{\ast }\right) ^{-} \\
& =\bigcap_{H\in \mathcal{H}}\bigcap_{\mu \in \mathbb{R}_{+}^{H}}\left[
\left( f+\sum_{t\in H}\mu _{t}f_{t}\right) _{\infty }\leq 0\right] \text{
(by (\ref{marco2})} \\
& =\bigcap_{H\in \mathcal{H}}\bigcap_{\mu \in \mathbb{R}_{+}^{H}}\left[
\left( f_{\infty }+\sum_{t\in H}\mu _{t}(f_{t})_{\infty }\right) \leq 0%
\right] \text{ (by (\ref{marco4})} \\
& =\bigcap_{H\in \mathcal{H}}\left[ \left( \sup_{\mu \in \mathbb{R}%
_{+}^{H}}\left( f_{\infty }+\sum_{t\in H}\mu _{t}(f_{t})_{\infty }\right)
\right) \leq 0\right] \\
& =\bigcap_{H\in \mathcal{H}}\left[ \left( f_{\infty }+\sup_{\mu \in \mathbb{%
R}_{+}^{H}}\sum_{t\in H}\mu _{t}(f_{t})_{\infty }\right) \leq 0\right] \\
& =\bigcap_{H\in \mathcal{H}}\left[ \left( f_{\infty }+\delta _{\left[
\sup_{t\in H}(f_{t})_{\infty }\leq 0\right] }\right) \leq 0\right]
=\bigcap_{H\in \mathcal{H}}\bigcap_{t\in H}\left[ (f_{t})_{\infty }\leq 0%
\right] \cap \lbrack f_{\infty }\leq 0] \\
& =\bigcap_{t\in T}\left[ (f_{t})_{\infty }\leq 0\right] \cap \lbrack
f_{\infty }\leq 0]=(\mathrm{P})_{\infty },
\end{align*}%
the penultimate equality coming from the fact that the family $\mathcal{H}$
is covering. We conclude the proof of Theorem \ref{thm4} with Lemmas \ref%
{lem3} and \ref{lem6}.\hfill $\square $

\begin{remark}
\label{rem6} Note that if $X=X^{\ast }=\mathbb{R}^{n}$ then the function $%
\varphi _{\mathcal{H}}$, when convex, is automatically $\tau (X^{\ast },X)$%
-quasicontinuous since any extended real-valued convex function on $\mathbb{R%
}^{n}$ with non-empty domain is quasicontinuous (e.g., \cite[Theorem 10.1]%
{Rocka1970}).
\end{remark}

\begin{corollary}[$\mathcal{H}_{1}$-reverse strong duality]
\label{cor5} Assume that $(\mathrm{P})$ satisfies the following conditions:%
\newline
$\mathrm{(a)}$\textrm{\ }$\func{dom}f\subset \bigcap\nolimits_{t\in T}\func{%
dom}f_{t},$\newline
$\mathrm{(b)}$\textrm{\ }$T$ is a convex and compact subset of some locally
convex topological vector space, \newline
$\mathrm{(c)}$\textrm{\ }$T\ni t\mapsto f_{t}(x)$ is concave and usc on $T$
for each $x\in \bigcap\nolimits_{t\in T}\func{dom}f_{t},$\newline
$\mathrm{(d)}\ $There exists$\ (\bar{t},\bar{\mu})\in T\times \mathbb{R}%
_{+}\,$\ such that $f+\bar{\mu}f_{\bar{t}}$ is weakly inf-locally compact,
\newline
$\mathrm{(e)}$\ $(\mathrm{P})_{\infty }$ is a linear subspace, \newline
Then,
\begin{equation*}
\min (\mathrm{P})=\sup_{(t,\mu )\in T\times \mathbb{R}_{+}}\inf_{x\in X}%
\big\{f(x)+\mu f_{t}(x)\big\}\in \mathbb{R}.
\end{equation*}
\end{corollary}

\noindent \textbf{Proof }From the first three assumptions and \cite[Remark
5.5]{DGLV21} we get that $\varphi _{\mathcal{H}_{1}}$ is convex. Moreover, $%
\varphi _{\mathcal{H}_{1}}=\inf_{(t,\mu )\in T\times \mathbb{R}_{+}}(f+\mu
f_{t})^{\ast }$ is majorized by the function $(f+\bar{\mu}f_{\bar{t}})^{\ast
}$, which is $\tau (X^{\ast },X)$-quasicontinuous by Lemma \ref{lem4} as, by
$\mathrm{(d)}$, $f+\bar{\mu}f_{\bar{t}}\in \Gamma (X)$ is weakly inf-locally
compact. So, by Lemma \ref{lem5}, $\varphi _{\mathcal{H}_{1}}$ is $\tau
(X^{\ast },X)$-quasicontinuous, and we conclude the proof by applying
Theorem \ref{thm4} with $\mathcal{H}=\mathcal{H}_{1}$ thanks to $\mathrm{(e).%
}$ \hfill $\square $

The next result recovers a variant of the reverse duality theorem of \cite[%
Theorem 3.3]{GLV14}{.}

\begin{corollary}[$\mathcal{F}(T)$-reverse strong duality]
\label{cor6} Assume that $E\cap \func{dom}f\neq \emptyset $ and that the two
following conditions are satisfied:%
\begin{equation*}
\begin{tabular}{l}
$\mathrm{(f)}\ \ \exists \lambda \in \mathbb{R}_{+}^{(T)}\text{ such that }%
\,f+\sum_{t\in T}\lambda _{t}f_{t}\text{ is weakly inf-locally compact.}$ \\
$\mathrm{(e)}\ \ (\mathrm{P})_{\infty }\text{ is a linear subspace.}$%
\end{tabular}%
\end{equation*}%
\newline
Then we have%
\begin{equation*}
\min (\mathrm{P})=\sup (\mathrm{D})=\sup_{\lambda \in \mathbb{R}%
_{+}^{(T)}}\inf_{x\in X}\left\{ f(x)+\sum_{t\in T}\lambda
_{t}f_{t}(x)\right\} \in \mathbb{R}.
\end{equation*}
\end{corollary}

\noindent \textbf{Proof }Condition $\mathrm{(f)}$ amounts to
\begin{equation*}
\exists H\in \mathcal{F}(T),\exists \mu \in \mathbb{R}_{+}^{H}\text{ such
that }\,f+\sum_{t\in H}\mu _{t}f_{t}\text{ weakly inf-locally compact.}
\end{equation*}%
Moreover, $\varphi _{\mathcal{F}(T)}$ is majorized by $\left( f+\sum_{t\in
H}\mu _{t}f_{t}\right) ^{\ast }$ which is $\tau (X^{\ast },X)$%
-quasicontinuous by Lemma \ref{lem4}. By Lemma \ref{lem5}, $\varphi _{%
\mathcal{F}(T)}$ is then $\tau (X^{\ast },X)$-quasicontinuous. Taking $%
\mathcal{H}=\mathcal{F}(T)$ in Theorem \ref{thm4} we obtain, by (\ref{2.1})
and (\ref{2.2}),
\begin{equation*}
\min (\mathrm{P})=\sup (\mathrm{D}_{\mathcal{H}})=\sup (\mathrm{D}),
\end{equation*}%
and the proof is complete. \hfill $\square $

We finally consider the countable case when $T=\mathbb{N}.$ Let $(\mathrm{P}%
_{\mathbb{N}}),$ $(\mathrm{P}_{m}),$ $(\mathrm{D}_{\mathbb{N}}),$ and $(%
\mathrm{D}_{m})$ be as in (\ref{2.7}), (\ref{2.8}), (\ref{2.9}), and (\ref%
{2.10}), respectively.

\begin{corollary}[$\mathcal{H}_{\mathbb{N}}$-reverse strong duality]
\label{cor7} Assume $\inf (\mathrm{P}_{\mathbb{N}})\neq +\infty $ and the
two conditions below are satisfied:%
\begin{equation*}
\begin{tabular}{l}
$\mathrm{(g)}\ \ \exists \ (N,\mu )\in \mathbb{N}\times \mathbb{R}_{+}^{N}%
\text{ such that }\,f+\sum_{k=1}^{N}\mu _{k}f_{k}\text{ is weakly
inf-locally compact,}$ \\
$\mathrm{(e)}\ \ (\mathrm{P})_{\infty }\text{ is a linear subspace.}$%
\end{tabular}%
\end{equation*}%
Then we have
\begin{equation*}
\min (\mathrm{P}_{\mathbb{N}})=\lim_{m\rightarrow \infty }\inf (\mathrm{P}%
_{m})=\lim_{m\rightarrow \infty }\sup (\mathrm{D}_{m})=\sup (\mathrm{D}_{%
\mathbb{N}}).
\end{equation*}%
Moreover, the optimal solution set of $(\mathrm{P}_{\mathbb{N}})$ is the sum
of a weakly compact convex set and a finite dimensional linear subspace.
\end{corollary}

\noindent \textbf{Proof} Since the covering family $\mathcal{H}_{\mathbb{N}}$
is directed we know that $\varphi _{\mathcal{H}_{\mathbb{N}}}$ is a convex
function. Moreover, $\varphi _{\mathcal{H}_{\mathbb{N}}}$ is majorized by $%
\left( f+\sum_{k=1}^{N}\mu _{k}f_{k}\right) ^{\ast }$ which is $\tau
(X^{\ast },X)$-quasicontinuous by Lemma \ref{lem4}. By Lemma \ref{lem5}, $%
\varphi _{\mathcal{H}_{\mathbb{N}}}$ is then $\tau (X^{\ast },X)$%
-quasicontinuous and, {by }\cite[Formula (5.6)]{DGLV21}, $\sup (\mathrm{D}_{%
\mathbb{N}})=\lim_{m\rightarrow \infty }\sup (\mathrm{D}_{m}).$ Applying
Theorem \ref{thm4} with $\mathcal{H}=\mathcal{H}_{\mathbb{N}}$ we obtain,
\begin{equation*}
\min (\mathrm{P}_{\mathbb{N}})=\sup (\mathrm{D}_{\mathbb{N}})=\sup_{m\in
\mathbb{N}}\sup (\mathrm{D}_{m})=\lim_{m\rightarrow +\infty }\sup (\mathrm{D}%
_{m})\leq \lim_{m\rightarrow +\infty }\inf (\mathrm{P}_{m})\leq \min (%
\mathrm{P}_{\mathbb{N}}),
\end{equation*}%
and the proof is complete. \hfill $\square $

\begin{remark}
\label{rem5} We now comment conditions $\mathrm{(a)-(g)}$ when $X=\mathbb{R}%
^{n},$ that is, in CSIP. Conditions $\mathrm{(d),}$ $\mathrm{(f)}$, and $%
\mathrm{(g)}$ are obviously satisfied while condition $\mathrm{(e)}$ is
equivalent \cite[Exercise 8.15]{GL98} to\newline
\indent$\mathrm{(h)}$ $f_{\infty }\left( x\right) >0,\forall x\in \left[
\left( 0^{+}E\right) \cap M^{\bot }\right] \diagdown \left\{ 0_{n}\right\} ,$%
\newline
where $M=\left\{ x\in \limfunc{lin}\left( 0^{+}E\right) :f_{\infty }\left(
x\right) =0=f_{\infty }\left( -x\right) \right\} .$ So, Corollary \ref{cor6}
is, in the CSIP setting, equivalent to {\cite[Theorem 3.2]{Karney1983} (see
also \cite[Theorem 8.8(i)]{GL98})}. Analogously,\ {\cite[Corollary 4.2]%
{Karney1983}} is the CSIP version of Corollary \ref{cor7}. \newline
If $(\mathrm{P})$ is the LSIP problem in (\ref{1.4}), we can write $%
f(x)=\left\langle c^{\ast },x\right\rangle $ and$\ f_{t}(x)=\left\langle
a_{t}^{\ast },x\right\rangle -b_{t},\;t\in T.$ Then since all functions have
full domain, $\mathrm{(a)}$\ trivially holds. Moreover, since $(\mathrm{P}%
)_{\infty }=\bigcap_{t\in T}[a_{t}^{\ast }\leq 0]\cap \lbrack c^{\ast }\leq
0],$ condition $\mathrm{(e)}$ can be expressed as follows:\newline
\indent$\mathrm{(e}^{\prime }\mathrm{)}$ $\left\{ x\in X:\left\langle
c^{\ast },x\right\rangle \leq 0;\ \left\langle a_{t}^{\ast },x\right\rangle
\leq 0,\forall t\in T\right\} $ is a linear subspace.\newline
Taking into account that a convex cone $K$ is a subspace if and only if $%
-K\subset K,$ $\mathrm{(e}^{\prime }\mathrm{)}$ is equivalent to\newline
\indent$\mathrm{(e}^{\prime \prime }\mathrm{)}$ $\left[ \left\langle c^{\ast
},x\right\rangle \leq 0;\ \left\langle a_{t}^{\ast },x\right\rangle \leq
0,\forall t\in T\right] \Longrightarrow \left[ \left\langle c^{\ast
},x\right\rangle =0=\ \left\langle a_{t}^{\ast },x\right\rangle ,\forall
t\in T\right] .$\newline
Moreover, condition $\mathrm{(e}^{\prime }\mathrm{)}$ can be reformulated in
terms of the data as \newline
\indent$\mathrm{(e}^{\prime \prime \prime }\mathrm{)}$ The pointed cone of $%
\overline{\limfunc{cone}}\left( \left\{ c^{\ast };a_{t}^{\ast },t\in
T\right\} \times \mathbb{R}_{+}\right) $ (i.e., its intersection with the
orthogonal subspace to its lineality) is a half-line in $\mathbb{R}^{n+1}$
\cite[Theorem 5.13(ii)]{GL98} (or, more precisely, the half-line $\mathbb{R}%
_{+}\left( 0_{n},1\right) $ \cite[page 155]{GJR02}).\newline
In the same vein, since $\func{dom}f=\mathbb{R}^{n},$ $f_{\infty
}=\left\langle c^{\ast },\cdot \right\rangle ,$ $0^{+}E=\bigcap_{t\in
T}[a_{t}^{\ast }\leq 0],$ and
\begin{equation*}
M^{\bot }=\left\{ x:\left\langle c^{\ast },x\right\rangle =0=\left\langle
a_{t}^{\ast },x\right\rangle ,\forall t\in T\right\} ^{\bot }=\limfunc{span}%
\left\{ c^{\ast };a_{t}^{\ast },t\in T\right\} ,
\end{equation*}%
condition $\mathrm{(h)}$ can be expressed as\newline
\indent$\mathrm{(h}^{\prime }\mathrm{)}$ $\left\langle c^{\ast
},x\right\rangle >0,\forall x\in \left( \bigcap_{t\in T}[a_{t}^{\ast }\leq
0]\right) \cap \limfunc{span}\left\{ c^{\ast };a_{t}^{\ast },t\in T\right\}
\diagdown \left\{ 0_{n}\right\} .$\newline
\end{remark}

\begin{example}
\label{Exam1}Consider the linear semi-infinite programming problem
\begin{equation*}
\begin{tabular}{lll}
$\left( \mathrm{P}\right) $ & $\inf\limits_{x\in \mathbb{R}^{2}}$ & $f\left(
x\right) =\left\langle c^{\ast },x\right\rangle $ \\
& s.t. & $-tx_{1}+(t-1)x_{2}+t-t^{2}\leq 0,\;t\in \left[ 0,1\right] ,$%
\end{tabular}%
\end{equation*}%
with $c^{\ast }\in \mathbb{R}_{+}^{2}\diagdown \left\{ \left( 0,0\right)
\right\} $ (see \cite[Example 3.1]{DGLV21}). According to Remark \ref{rem5},
$\mathrm{(a),}$ $\mathrm{(d),}$ $\mathrm{(f)}$, and $\mathrm{(g)}$ hold
independently of the data. Condition $\mathrm{(b)}$ holds because $\left[ 0,1%
\right] \subset \mathbb{R}$ is compact and convex and $\mathrm{(c)}$ because
$t\longmapsto -tx_{1}+(t-1)x_{2}+t-t^{2}$\ is concave on $\mathbb{R}$\ for
any $x\in \mathbb{R}^{2}.$ Regarding $\mathrm{(e),}$ the set in $\mathrm{(e}%
^{\prime }\mathrm{)}$
\begin{equation*}
\left\{ x\in \mathbb{R}^{2}:\left\langle c^{\ast },x\right\rangle \leq
0;-tx_{1}+(t-1)x_{2}\leq 0,\forall t\in \left[ 0,1\right] \right\} =\left\{
x\in \mathbb{R}_{+}^{2}:\left\langle c^{\ast },x\right\rangle \leq 0\right\}
\end{equation*}%
is $\left\{ \left( 0,0\right) \right\} $ when $c^{\ast }$ belongs to the
interior $\mathbb{R}_{++}^{2}$ of $\mathbb{R}_{+}^{2}$ and a positive axis
when $c^{\ast }$ belongs to its boundary. Hence, $\mathrm{(e)}$ only holds
for $c^{\ast }\in \mathbb{R}_{++}^{2}.$ Observe that the cone in $\mathrm{(e}%
^{\prime \prime }\mathrm{)}$ is
\begin{equation*}
\limfunc{cone}\left\{ \left(
\begin{array}{c}
c_{1}^{\ast } \\
c_{2}^{\ast }%
\end{array}%
\right) ,\left(
\begin{array}{c}
-1 \\
0%
\end{array}%
\right) ,\left(
\begin{array}{c}
0 \\
-1%
\end{array}%
\right) \right\} \times \mathbb{R}_{+},
\end{equation*}%
and its pointed cone is
\begin{equation*}
\mathbb{R}_{+}\left(
\begin{array}{c}
0 \\
0 \\
1%
\end{array}%
\right) \text{ }\left( \text{resp., }\limfunc{cone}\left\{ \left(
\begin{array}{c}
-1 \\
0 \\
0%
\end{array}%
\right) ,\left(
\begin{array}{c}
0 \\
0 \\
1%
\end{array}%
\right) \right\} ,\limfunc{cone}\left\{ \left(
\begin{array}{c}
0 \\
-1 \\
0%
\end{array}%
\right) ,\left(
\begin{array}{c}
0 \\
0 \\
1%
\end{array}%
\right) \right\} \right) ,
\end{equation*}%
when $c^{\ast }\in \mathbb{R}_{++}^{2}$ ($c^{\ast }\in \mathbb{R}_{++}\left(
1,0\right) ,c^{\ast }\in \mathbb{R}_{++}\left( 1,0\right) ,$ respectively).
So, we get again that $\mathrm{(e)}$ only holds for $c^{\ast }\in \mathbb{R}%
_{++}^{2}.$ Regarding condition $\mathrm{(h)}$, since $\bigcap_{t\in \left[
0,1\right] }[a_{t}^{\ast }\leq 0]=\mathbb{R}_{+}^{2},$ and $\limfunc{span}%
\left\{ c^{\ast };a_{t}^{\ast },t\in T\right\} =\mathbb{R}^{2}$ if $c^{\ast
}\in \mathbb{R}_{++}^{2},$ so that $\mathrm{(h)}$ holds while $\limfunc{span}%
\left\{ c^{\ast };a_{t}^{\ast },t\in T\right\} $ is a positive axis and $%
\mathrm{(h)}$ fails, otherwise. Thus, $\mathrm{(e)}$ and $\mathrm{(h)}$ hold
or not simultaneously. \newline
In conclusion, by Corollary \ref{cor5}, $\mathcal{H}_{1}$-reverse strong
duality holds whenever $c^{\ast }\in \mathbb{R}_{++}^{2}$ while, by
Corollary \ref{cor6}, $\mathcal{F}(T)$-reverse strong duality holds whenever
$c^{\ast }\in \mathbb{R}_{++}^{2}.$ Observe that, from the direct
computations carried out in \cite[Example 3.1]{DGLV21}, $\mathcal{H}_{1}$%
-reverse strong duality actually holds for all $c^{\ast }\in \mathbb{R}%
_{+}^{2}\diagdown \left\{ \left( 0,0\right) \right\} .$
\end{example}

\begin{example}
\label{Exam2}The countable linear semi-infinite programming problem%
\begin{equation*}
\begin{array}{lll}
\left( \mathrm{P}_{\mathbb{N}}\right) & \inf\limits_{x\in \mathbb{R}^{2}} &
x_{2}\,\  \\
& \text{\textit{s.t.}} & x_{1}+k\left( k+1\right) x_{2}\geq 2k+1,\;k\in
\mathbb{N},%
\end{array}%
\end{equation*}%
violates the assumptions of Corollaries \ref{cor5}, \ref{cor6}, and \ref%
{cor7}, as $\mathrm{(b)}$ and $\mathrm{(c)}$\ obviously fail, as well as $%
\mathrm{(e)}$ and $\mathrm{(h)}$. In fact, $\mathrm{(e}^{\prime }\mathrm{)}$
and $\mathrm{(e}^{\prime \prime }\mathrm{)}$ fail because
\begin{equation*}
\left\{ x\in \mathbb{R}^{2}:x_{2}\leq 0,-x_{1}-k\left( k+1\right) x_{2}\leq
0,\;k\in \mathbb{N}\right\} =\mathbb{R}_{+}\times \left\{ 0\right\}
\end{equation*}%
is not a linear subspace and the pointed cone of%
\begin{equation*}
\overline{\limfunc{cone}}\left\{ \left( 0,1\right) ;\left( -1,-k\left(
k+1\right) \right) ,k\in \mathbb{N}\right\} \times \mathbb{R}_{+}=\left\{
x\in \mathbb{R}^{3}:x_{1}\leq 0,x_{3}\geq 0\right\}
\end{equation*}%
is not a half-line, respectively, while $\mathrm{(h)}$ fails because $x_{2}$%
~vanishes on an edge of
\begin{equation*}
\left( 0^{+}E\right) \cap M^{\bot }=0^{+}E\cap \mathbb{R}^{2}=\limfunc{cone}%
\left\{ \left( -2,1\right) ,\left( 1,0\right) \right\} \newline
.
\end{equation*}%
So, we cannot apply the mentioned corollaries to conclude that $\mathcal{H}$%
-reverse strong duality holds for $\mathcal{H=\mathcal{H}}_{1},\mathcal{H}_{%
\mathbb{N}},\mathcal{F}(T).$ Actually, $\mathcal{H}$-reverse strong duality
does not hold for these three families because the feasible set of $\left(
\mathrm{P}_{\mathbb{N}}\right) $ is
\begin{equation*}
E=\limfunc{co}\left( \left\{ \left( k,\frac{1}{k}\right) ,k\in \mathbb{N}%
\right\} \cup \left\{ x\in \mathbb{R}^{2}:x_{1}+2x_{2}=3,x_{1}\leq 1\right\}
\right) ,
\end{equation*}%
which implies $\inf (\mathrm{P}_{\mathbb{N}})=0$ with $\limfunc{sol}(\mathrm{%
P}_{\mathbb{N}})=\emptyset ,$ while $\sup \left( D\right) =-\infty ,$ which
in turn implies $\sup \left( D_{\mathcal{H}}\right) =-\infty $ for any $%
\mathcal{H}$ such that $\emptyset \neq \mathcal{H\subset F}(T),$ by (\ref%
{2.5}).
\end{example}

\section{ $\mathcal{H}\mathcal{-}$Farkas lemma}

We now establish some new versions of Farkas lemma relative to a given
family $\mathcal{H}\subset \mathcal{F}(T)$. These results assert the
equivalence between some inclusion $\left( \mathrm{i}\right) $ of the
solution set $E$ of $\left\{ f_{t}(x)\leq 0,t\in T\right\} $ into certain
set involving $f$ and some condition $\left( \mathrm{ii}\right) $ involving $%
\{f;\ f_{t},t\in T\}$ and $\mathcal{H}.$ We first provide a Farkas-type
result \ relative to the family $\mathcal{H}_{1}$ without assuming the lower
semicontinuity of the involved functions. Stronger results
(characterizations of Farkas lemma) will be then obtained under the lower
semicontinuity (or even continuity) assumption.

\begin{proposition}[$\mathcal{H}_{1}$-Farkas lemma]
\label{lem1} Assume that conditions $\mathrm{(a)}$,$\mathrm{(b)}$,$\mathrm{%
(c)}$ in Corollary \ref{cor5} altogether with the generalized Slater
condition:
\begin{equation*}
\exists \bar{x}\in \func{dom}f:\ \ f_{t}(\bar{x})<0,\ \forall t\in T.
\end{equation*}

Then, for any $\alpha \in \mathbb{R}$, the following statements are
equivalent:\newline
$\left( \mathrm{i}\right) $\ \ $\left[ f_{t}(x)\leq 0,\forall t\in T\right]
\ \Longrightarrow \ f(x)\geq \alpha .$\newline
$\left( \mathrm{ii}\right) $ There exist $\bar{t}\in T$ and $\bar{\mu}\in
\mathbb{R}_{+}$ such that
\begin{equation}
f(x)+\bar{\mu}f_{\bar{t}}(x)\geq \alpha , \ \ \forall x\in X.
\label{marco321}
\end{equation}
\end{proposition}

\noindent \textbf{Proof } We observe first that $\left( \mathrm{i}\right) $
is equivalent to $\inf (\mathrm{P})\geq \alpha ,$ where $(\mathrm{P})$ is
the CIP in (\ref{1.1}). So, it follows from \cite[Theorem 5.3]{DGLV21} that $%
\inf (\mathrm{P})=\max (\mathrm{D}_{\mathcal{H}_{1}})\geq \alpha $; i.e., $%
\left( \mathrm{i}\right) $ is equivalent to
\begin{equation*}
\max\limits_{(t,\mu )\in T\times \mathbb{R}}\inf_{x\in \func{dom}f}\{f(x)+{%
\mu }f_{t}(x)\}\geq \alpha .
\end{equation*}%
In other words, there exists $(\bar{t},\bar{\mu})\in T\times \mathbb{R}$
satisfying (\ref{marco321}), which is $\left( \mathrm{ii}\right) $, and we
are done. \hfill \hfill \noindent $\square \medskip $

Observe that statement $\left( \mathrm{i}\right) $ means that $E$ is
contained in the reverse convex set $\left\{ x\in X:f(x)\geq \alpha \right\}
$ while $\left( \mathrm{ii}\right) $ would be the same replacing the
infinite family $\left\{ f_{t},t\in T\right\} $ by the singleton one $%
\left\{ f_{\bar{t}}\right\} ,$ so that Lemma \ref{lem1} characterizes when
an inequality $f(x)\geq \alpha $ is consequence of some single constraint $%
f_{t}(x)\leq 0.$

The next two propositions provide, under the lower semicontinuity
assumption, a characterization in terms of $\mathcal{A}_{\mathcal{H}}$
(statement $\left( \mathrm{I}\right) $) of the Farkas lemma (statement $%
\left( \mathrm{II}\right) $) relative to an arbitrary non-empty covering
family $\mathcal{H}\subset \mathcal{F}(T)$.

\begin{proposition}[Characterization of $\mathcal{H}$-Farkas lemma]
\label{lem2} Let $\mathcal{H}\subset \mathcal{F}(T)$ be a covering family.
Assume that $\{f;\ f_{t},t\in T\}\subset \Gamma (X)$, $E\cap (\func{dom}%
f)\neq \emptyset $, and consider the following statements:

$\left( \mathrm{I}\right) $ $\mathcal{A}_{\mathcal{H}}$ is $w^{\ast }$%
-closed convex regarding $\{0_{X^{\ast }}\}\times \mathbb{R}$.

$\left( \mathrm{II}\right) $ For $\alpha \in \mathbb{R}$, the next two
conditions are equivalent:

$\qquad \left( \mathrm{i}\right) $ $\left[ f_{t}(x)\leq 0,\forall t\in T%
\right] \Longrightarrow \ f(x)\geq \alpha ,$

$\qquad \left( \mathrm{ii}\right) $ there exist $H\in \mathcal{H}\ $and $\mu
\in \mathbb{R}_{+}^{H}$ such that
\begin{equation}
f(x)+\sum\limits_{t\in H}\mu _{t}f_{t}(x)\geq \alpha ,\forall x\in X.
\label{marco333}
\end{equation}%
Then, $[(\mathrm{I})\Longrightarrow (\mathrm{II})]$, and the converse
implication, $[(\mathrm{II})\Longrightarrow (\mathrm{I})],$ holds when%
\textbf{\ }$\inf (\mathrm{P})\in \mathbb{R}$.
\end{proposition}

\noindent \textbf{Proof } By the characterization of $\mathcal{H}$-strong
duality at a point in (\ref{4.4}), applied to $x^{\ast }=0_{X^{\ast }}$, one
gets that $\left( \mathrm{I}\right) $ is equivalent to
\begin{equation}  \label{eqnewa}
\inf (\mathrm{P})=\max (\mathrm{D}_{\mathcal{H}}),
\end{equation}%
which is itself equivalent to the existence of $H\in \mathcal{H}\ $and $\mu
\in \mathbb{R}_{+}^{H}$ such that
\begin{equation*}
\inf (\mathrm{P})=\inf_{x\in X}\left( f(x)+\sum\limits_{t\in H}\mu
_{t}f_{t}(x\right) .
\end{equation*}%
Since $\left( \mathrm{i}\right) $ is equivalent to $\inf (\mathrm{P})\geq
\alpha ,$ it now follows that $[(\mathrm{I})\Longrightarrow (\mathrm{II})]$.

Conversely, if $\inf (\mathrm{P})\in \mathbb{R}$ and $(\mathrm{II})$ holds,
then just take $\alpha =\inf (\mathrm{P})$. As $(\mathrm{II})$ holds, it
follows that there are $H\in \mathcal{H}$ and $\mu \in \mathbb{R}_{+}^{H}$
such that (\ref{marco333}) holds, and
\begin{equation*}
\sup (\mathrm{D}_{\mathcal{H}})\geq \inf_{x\in X}\left(
f(x)+\sum\limits_{t\in H}\mu _{t}f_{t}(x)\right) \geq \alpha =\inf (\mathrm{P%
}).
\end{equation*}%
In other words, $\sup (\mathrm{D}_{\mathcal{H}})=\inf (\mathrm{P}),$ $\sup (%
\mathrm{D}_{\mathcal{H}})$ is attained at $H\in \mathcal{H}\ $and $\mu \in
\mathbb{R}_{+}^{H},$ meaning that \eqref{eqnewa} holds, which is  $(\mathrm{I%
}) $, and the proof is complete. $\ \hfill \square $

\begin{remark}
In the special case when $\mathcal{H}=\mathcal{F}(T)$ the condition $(%
\mathrm{ii})$ in Proposition \ref{lem2} reads as

$\qquad (\mathrm{ii}^{\prime })$ there exists $\lambda \in \mathbb{R}%
_{+}^{(T)}$ such that $f(x)+\sum\limits_{t\in T}\lambda _{t}f_{t}(x)\geq
\alpha ,$ for all $x\in X$,\newline
and Proposition \ref{lem2} goes back to the Farkas lemma given in \cite[%
Theorem 2]{DGLS07} under a slightly different qualification condition. So,
Proposition \ref{lem2} is a variant of \cite[Theorem 2]{DGLS07}.
\end{remark}

Let us get back to the linear case, where
\begin{equation}
f(x)=\left\langle c^{\ast },x\right\rangle ,\ f_{t}(x)=\left\langle
a_{t}^{\ast },x\right\rangle -b_{t},t\in T,  \label{4.2}
\end{equation}%
with $\left\{ c^{\ast };a_{t}^{\ast },t\in T\right\} \subset X^{\ast },$ and
$\left\{ b_{t},t\in T\right\} \subset \mathbb{R}.$ Then, $\mathcal{A}_{%
\mathcal{H}}=\left\{ \left( c^{\ast },0\right) \right\} +\mathcal{K}_{%
\mathcal{H}}$ (see \cite[(4.4)]{DGLV21}), where

\begin{equation*}
\mathcal{K}_{\mathcal{H}}=\tbigcup\limits_{H\in \mathcal{H}}\limfunc{cone}%
\left( \left\{ (a_{t}^{\ast },b_{t}),\ t\in H\right\} +\left\{ 0_{X^{\ast
}}\right\} \times \mathbb{R}_{+}\right) .
\end{equation*}%
In particular,
\begin{equation*}
\begin{array}{ll}
\mathcal{K}_{\mathcal{H}_{1}} & =\tbigcup\limits_{t\in T}\limfunc{cone}%
\left\{ (a_{t}^{\ast },b_{t}+\varepsilon ):\varepsilon \geq 0\right\}%
\end{array}%
\end{equation*}%
and, by \cite[Proposition 4.1]{DGLV21},
\begin{equation*}
\mathcal{K}_{\mathcal{F}(T)}=\limfunc{cone}\left( \left\{ (a_{t}^{\ast
},b_{t}),\ t\in T\right\} +\left\{ 0_{X^{\ast }}\right\} \times \mathbb{R}%
_{+}\right) .
\end{equation*}%
For instance, for the LSIP problem\ in Example \ref{Exam1},%
\begin{equation*}
\begin{array}{ll}
\mathcal{K}_{\mathcal{H}_{1}} & =\tbigcup\limits_{t\in \left[ 0,1\right] }%
\limfunc{cone}\left\{ \left( -t,t-1,t^{2}-t+\varepsilon \right) :\varepsilon
\geq 0\right\}%
\end{array}%
\end{equation*}%
while $\mathcal{K}_{\mathcal{F}(T)}$ is (see \cite[Example 4.1]{DGLV21}) the
union of the origin with the epigraph of the convex function
\begin{equation*}
\psi \left( x\right) :=\left\{
\begin{array}{ll}
\frac{x_{1}x_{2}}{x_{1}+x_{2}}, & x\in \mathbb{R}_{-}^{2}\diagdown \left\{
0_{2}\right\} , \\
+\infty , & \text{else.}%
\end{array}%
\right.
\end{equation*}

We finish this section with a characterization, in terms of $\mathcal{K}_{%
\mathcal{H}},$ of the Farkas lemma (statement $\left( \mathrm{II}\right) $
below) relative to an arbitrary non-empty covering family $\mathcal{H}%
\subset \mathcal{F}(T)$.

\begin{proposition}[$\mathcal{H}$-Farkas lemma for linear infinite systems]
\label{lem73} Consider the linear functions $\{f;$~$f_{t},t\in T\}$ defined
in \eqref{4.2}, and suppose that $\inf (\mathrm{P})$ is finite and that $%
\mathcal{H}$ is a covering family. Given $c^{\ast }\in X^{\ast }$, the
following statements are equivalent:

$(\mathrm{I})$ $\overline{\limfunc{co}}(\mathcal{K}_{\mathcal{H}})\cap
(\{-c^{\ast }\}\times \mathbb{R}_{+})=\mathcal{K}_{\mathcal{H}}\cap
(\{-c^{\ast }\}\times \mathbb{R}_{+}),$

$(\mathrm{II})$ For $\alpha \in \mathbb{R}$, the following statements are
equivalent:

$\qquad (\mathrm{i})$ $\left[ \langle a_{t}^{\ast },x\rangle \leq
b_{t},\forall t\in T\right] \Longrightarrow \ \langle c^{\ast },x\rangle
\geq \alpha .$

$\qquad (\mathrm{ii})$ There exist $H\in \mathcal{H}$ and $\mu \in \mathbb{R}%
_{+}^{H}$ such that $\sum\limits_{t\in H}\mu _{t}a_{t}^{\ast }=-c^{\ast }$
and\ \ $-\sum\limits_{t\in H}\mu _{t}b_{t}\geq \alpha .$
\end{proposition}

\noindent \textbf{Proof.} When $\mathcal{H}$ is a covering family and $%
E\not=\emptyset $,\ according to \cite[Corollary 5.3]{DGLV21}, one has
\begin{equation}
\Big(\inf (\mathrm{P})=\max (\mathrm{D}_{\mathcal{H}}) \Big) %
\Longleftrightarrow \Big( \left( \overline{\limfunc{co}}\ \mathcal{K}_{%
\mathcal{H}}\right) \cap (\{-c^{\ast }\}\times \mathbb{R}_{+})=\mathcal{K}_{%
\mathcal{H}}\cap (\{-c^{\ast }\}\times \mathbb{R}_{+})\Big).  \label{4.3}
\end{equation}%
The rest of the proof is similar to that of Proposition \ref{lem2}, using (%
\ref{4.4}) and (\ref{4.3}). $\ \hfill \square \medskip $

\section{$\mathcal{H}$-optimality conditions}

In this section we establish optimality conditions for the problem $(\mathrm{%
P})$ associated with some family $\mathcal{H}\subset \mathcal{F}(T)$. We
shall represent by $\limfunc{sol}(\mathrm{D}_{\mathcal{H}})$ the set of
optimal solutions of\textbf{\ }$(\mathrm{D}_{\mathcal{H}}).$ In particular,
when $\mathcal{H}=\mathcal{F}(T),$ one obtains the classical KKT conditions
involving finitely many multipliers and, when $\mathcal{H}=\mathcal{H}_{1},$
optimality conditions involving a unique multiplier.

\begin{theorem}[Primal-dual $\mathcal{H}{-}$optimality condition]
\label{new1}Let $\bar{x}\in E\cap (\func{dom}f)$, $H\in \mathcal{H}$ and $%
\mu \in \mathbb{R}_{+}^{H}$. Then, the following statements are equivalent:%
\newline
$(\mathrm{i})$ $\bar{x}\in \limfunc{sol}(\mathrm{P}),$ $(H,\mu )\in \limfunc{%
sol}(\mathrm{D}_{\mathcal{H}}),$ and $\inf (\mathrm{P})=\sup (\mathrm{D}_{%
\mathcal{H}}).$\newline
$(\mathrm{ii})$ $f(\bar{x})=\inf_{X}\left( f+\tsum\limits_{t\in H}\mu
_{t}f_{t}\right) ,$ and $\mu _{t}f_{t}(\bar{x})=0,\ $for all $t\in H.$%
\newline
$(\mathrm{iii})$ $0_{X^{\ast }}\in \partial \left( f+\tsum\limits_{t\in
H}\mu _{t}f_{t}\right) (\overline{x}),\ $and $\mu _{t}f_{t}(\bar{x})=0,\ $%
for all $t\in H.$
\end{theorem}

\noindent \textbf{Proof. }$[(\mathrm{i})\Rightarrow (\mathrm{ii})]$ We have%
\begin{equation*}
\inf_{X}\left( f+\tsum\limits_{t\in H}\mu _{t}f_{t}\right) =\sup (\mathrm{D}%
_{\mathcal{H}})=\inf (\mathrm{P})=f(\bar{x}),
\end{equation*}%
and
\begin{equation*}
f(\bar{x})=\inf_{X}\left( f+\tsum\limits_{t\in H}\mu _{t}f_{t}\right) \leq f(%
\bar{x})+\tsum\limits_{t\in H}\mu _{t}f_{t}(\bar{x})\leq f(\bar{x}).
\end{equation*}%
Hence, $\tsum\limits_{t\in H}\mu _{t}f_{t}(\bar{x})=0$ and $(\mathrm{ii})$
holds.

$[(\mathrm{ii})\Rightarrow (\mathrm{iii})]$ We have%
\begin{equation*}
\left( f+\tsum\limits_{t\in H}\mu _{t}f_{t}\right) (\bar{x})=f(\bar{x}%
)=\inf_{X}\left( f+\tsum\limits_{t\in H}\mu _{t}f_{t}\right) .
\end{equation*}%
Thus, $\bar{x}\in \func{argmin}\left( f+\tsum\limits_{t\in H}\mu
_{t}f_{t}\right) $ or, equivalently, $0_{X^{\ast }}\in \partial \left(
f+\tsum\limits_{t\in H}\mu _{t}f_{t}\right) (\overline{x}).$

$[(\mathrm{iii})\Rightarrow (\mathrm{i})]$ Now we write%
\begin{equation*}
\inf (\mathrm{P})\leq f(\bar{x})=\left( f+\tsum\limits_{t\in H}\mu
_{t}f_{t}\right) (\bar{x})=\inf_{X}\left( f+\tsum\limits_{t\in H}\mu
_{t}f_{t}\right) \leq \sup (\mathrm{D}_{\mathcal{H}})\leq \inf (\mathrm{P}),
\end{equation*}%
and $(\mathrm{i})$ holds. $\hfill \square $

\begin{corollary}[1st $\mathcal{H}{-}$optimality condition for $(\mathrm{P})$%
]
\label{last}Assume that $\inf (\mathrm{P})=\max (\mathrm{D}_{\mathcal{H}})$
and let $\bar{x}\in E\cap (\func{dom}f)$. Then, the following statements are
equivalent:\newline
$(\mathrm{i})$ $\bar{x}\in \limfunc{sol}(\mathrm{P}).$\newline
$(\mathrm{ii})$ For each $(H,\mu )\in \limfunc{sol}(\mathrm{D}_{\mathcal{H}%
}),$ we have
\begin{equation}
0_{X^{\ast }}\in \partial \left( f+\tsum\limits_{t\in H}\mu _{t}f_{t}\right)
(\overline{x}),\ \text{and }\mu _{t}f_{t}(\bar{x})=0,\ \forall t\in H.
\label{marco322}
\end{equation}%
\newline
$(\mathrm{iii})$ There exists $(H,\mu )\in \limfunc{sol}(\mathrm{D}_{%
\mathcal{H}})$ such that (\ref{marco322}) is fulfilled.
\end{corollary}

\noindent \textbf{Proof. }$[(\mathrm{i})\Rightarrow (\mathrm{ii})]$ is just $%
[(\mathrm{i})\Rightarrow (\mathrm{iii})]$ in Theorem \ref{new1}.

$[(\mathrm{ii})\Rightarrow (\mathrm{iii})]$ is due to the assumption $%
\limfunc{sol}(\mathrm{D}_{\mathcal{H}})\neq \emptyset .$

$[(\mathrm{iii})\Rightarrow (\mathrm{i})]$ follows from $[(\mathrm{iii}%
)\Rightarrow (\mathrm{i})]$ in Theorem \ref{new1}. $\hfill \square $

\begin{corollary}[2nd $\mathcal{H}{-}$optimality condition for $(\mathrm{P})$%
]
\label{cor52} Let $\mathcal{H}\subset \mathcal{F}(T)$ be a covering family.
Assume that $\{f;\ f_{t},t\in T\}\subset \Gamma (X)$ and $E\cap (\func{dom}%
f)\neq \emptyset $. Assume further that $\mathcal{A}_{\mathcal{H}}$ is $%
w^{\ast }$-closed convex regarding $\{0_{X^{\ast }}\}\times \mathbb{R}$.
Then $\bar{x}\in $ $\limfunc{sol}(\mathrm{P})$ if and only if there exist $%
H\in \mathcal{H}$ and\ $\mu \in \mathbb{R}_{+}^{H}$ such that %
\eqref{marco322} holds.
\end{corollary}

\noindent \textbf{Proof } Taking $x^{\ast }=0_{X^{\ast }}$ in \eqref{4.4}\
one has $\inf \mathrm{(P)}=\max (\mathrm{D}_{\mathcal{H}}).$ Corollary \ref%
{last} concludes the proof. $\hfill \square $

\begin{remark}
When $\mathcal{H}=\mathcal{F}(T)$, the conclusion of Corollary \ref{cor52}
is that $\bar{x}\in \limfunc{sol}\left( \mathrm{P}\right) $ if and only if
there exist $\lambda \in \mathbb{R}_{+}^{(T)}$ such that
\begin{equation*}
0_{X^{\ast }}\in \partial \left( f+\sum\limits_{t\in T}\lambda
_{t}f_{t}\right) (\bar{x})\ \text{and}\ \lambda_{t}f_{t}(\bar{x})=0,\forall
t\in T,
\end{equation*}%
which recalls us about the optimality condition given in \cite[Theorem 3]%
{DGLS07} under the assumptions that both the sets $\mathcal{K}_{\mathcal{F}%
(T)}$ and $\limfunc{epi}f^{\ast }+\overline{\mathcal{K}_{\mathcal{F}(T)}}$
are $w^{\ast }$-closed.
\end{remark}

\begin{corollary}[$\mathcal{H}{-}$\textbf{optimality condition for linear }$(%
\mathrm{P})$]
\label{cor53} Let $(\mathrm{P})$ be linear with $E \not= \emptyset$. Let $%
\mathcal{H}$ be a covering family. Assume that $\mathcal{K}_{\mathcal{H}}$
is weak$^{\ast }$-closed convex regarding $\{-c^{\ast }\}\times \mathbb{R}$.
Then $\bar{x}\in \limfunc{sol}(\mathrm{P})$ if and only if there exist $H\in
\mathcal{H}$ and $\mu \in \mathbb{R}_{+}^{H}$ such that
\begin{equation}
\sum\limits_{t\in H}\mu _{t}a_{t}^{\ast }=-c^{\ast }\ \text{and}\ \
\sum\limits_{t\in H}\mu _{t}b_{t}=-\left\langle c^{\ast },\bar{x}%
\right\rangle .  \label{5.3}
\end{equation}
\end{corollary}

\noindent \textbf{Proof. } By \cite[Corollaty 5.3]{DGLV21} one has $\inf
\mathrm{(P)}=\max (\mathrm{D}_{\mathcal{H}})$. In the linear case one has %
\eqref{marco322} $\Leftrightarrow $ \eqref{5.3}. We conclude the proof with
Corollary \ref{last}. $\hfill \square $

\begin{corollary}[\textbf{Optimality condition for }$(\mathrm{D}_{\mathcal{H}%
})$]
\label{cor54}Assume that $\min (\mathrm{P})=\sup (\mathrm{D}_{\mathcal{H}%
})\neq +\infty $, and let $H\in \mathcal{H}$ and $\mu \in \mathbb{R}_{+}^{H}$%
. Then, the following statements are equivalent:\newline
$(\mathrm{i})$ $(H,\mu )\in \limfunc{sol}(\mathrm{D}_{\mathcal{H}}).$\newline
$(\mathrm{ii})$ For each $\bar{x}\in \limfunc{sol}(\mathrm{P}),$ (\ref%
{marco322}) holds.\newline
$(\mathrm{iii})$ There exists $\bar{x}\in \limfunc{sol}(\mathrm{P})$ such
that (\ref{marco322}) is fulfilled.
\end{corollary}

\noindent \textbf{Proof. }$[(\mathrm{i})\Rightarrow (\mathrm{ii})]$ follows
from $[(\mathrm{i})\Rightarrow (\mathrm{iii})]$ in Theorem \ref{new1}.

$[(\mathrm{ii})\Rightarrow (\mathrm{iii})]$ is due to the assumption $%
\limfunc{sol}(\mathrm{P})\neq \emptyset .$

$[(\mathrm{iii})\Rightarrow (\mathrm{i})]$ follows from $[(\mathrm{iii}%
)\Rightarrow (\mathrm{i})]$ in Theorem \ref{new1}. $\hfill \square \medskip $

We finish by revisiting again Example \ref{Exam1}, with $\mathcal{H=H}_{1}.$
For $c^{\ast }\in \mathbb{R}_{++}^{2},$ let us check the fulfilment of (\ref%
{5.3}) at $\bar{x}=\left( \left( \frac{c_{2}^{\ast }}{c_{1}^{\ast
}+c_{2}^{\ast }}\right) ^{2},\left( \frac{c_{1}^{\ast }}{c_{1}^{\ast
}+c_{2}^{\ast }}\right) ^{2}\right) .$ Taking $H=\left\{ \overline{t}%
\right\} ,$ with $\overline{t}=\frac{c_{1}^{\ast }}{c_{1}^{\ast
}+c_{2}^{\ast }}\in \left] 0,1\right[ ,$ and $\mu \in \mathbb{R}_{+}^{\left( %
\left[ 0,1\right] \right) }$ such that $\mu _{\overline{t}}=c_{1}^{\ast
}+c_{2}^{\ast }>0$ and $\mu _{t}=0$ for all $t\in \left[ 0,1\right]
\diagdown \left\{ \overline{t}\right\} ,$ one has%
\begin{equation*}
\sum\limits_{t\in H}\mu _{t}a_{t}^{\ast }=\left( c_{1}^{\ast }+c_{2}^{\ast
}\right) \left( -\frac{c_{1}^{\ast }}{c_{1}^{\ast }+c_{2}^{\ast }},-\frac{%
c_{2}^{\ast }}{c_{1}^{\ast }+c_{2}^{\ast }}\right) =-c^{\ast }
\end{equation*}%
and%
\begin{equation*}
\sum\limits_{t\in H}\mu _{t}b_{t}=\left( c_{1}^{\ast }+c_{2}^{\ast }\right)
\left( \left( \frac{c_{1}^{\ast }}{c_{1}^{\ast }+c_{2}^{\ast }}\right) ^{2}-%
\frac{c_{1}^{\ast }}{c_{1}^{\ast }+c_{2}^{\ast }}\right) =-\frac{c_{1}^{\ast
}c_{2}^{\ast }}{c_{1}^{\ast }+c_{2}^{\ast }}=-\left\langle c^{\ast },\bar{x}%
\right\rangle ,
\end{equation*}%
so that $\bar{x}\in \limfunc{sol}(\mathrm{P})$ (recall that $\mathcal{K}_{%
\mathcal{H}_{1}}$ is closed). Moreover, $(H,\mu )\in \limfunc{sol}(\mathrm{D}%
_{\mathcal{H}})$ by Corollary \ref{cor54} as
\begin{equation*}
\partial \left( c^{\ast }+\tsum\limits_{t\in H}\mu _{t}a_{t}^{\ast }\right)
=\left\{ c^{\ast }+\left( c_{1}^{\ast }+c_{2}^{\ast }\right) \left( -\frac{%
c_{1}^{\ast }}{c_{1}^{\ast }+c_{2}^{\ast }},-\frac{c_{2}^{\ast }}{%
c_{1}^{\ast }+c_{2}^{\ast }}\right) \right\} =\left\{ \left( 0,0\right)
\right\}
\end{equation*}%
and the complementarity condition $\mu _{t}f_{t}(\bar{x})=0,$ for all $t\in
T,$ holds.

\textbf{Acknowledgements}

This research was supported by Vietnam National University HoChiMinh city
(VNU-HCM) under the grant number B2021-28-03 (N. Dinh) and by Ministerio de
Ciencia, Innovaci\'{o}n y Universidades (MCIU), Agencia Estatal de
Investigaci\'{o}n (AEI), and European Regional Development Fund (ERDF),
Project PGC2018-097960-B-C22 (M.A. Goberna and M.A. L\'{o}pez).

\end{document}